\documentclass{amsart}[12pt]

\theoremstyle{plain}
\newtheorem{thm}{Theorem}[section]
\newtheorem{prop}[thm]{Proposition}
\newtheorem{lem}[thm]{Lemma}
\newtheorem{coro}[thm]{Corollary}

\newtheorem{ex}[thm]{Example}
\theoremstyle{definitions}
\newtheorem{dfn}[thm]{Definition}

\usepackage{amssymb}
\usepackage{tikz}
\usepackage{ytableau}

\begin{document}
\title[]{Descent Representations of Generalized Coinvariant Algebras}
\author[K. P. Meyer]{Kyle P. Meyer}
\address{Department of Mathematics, UCSD}
\email{kpmeyer@ucsd.edu}
\numberwithin{equation}{section}
\begin{abstract}
The coinvariant algebra $R_n$ is a well-studied $\mathfrak{S}_n$-module that is a graded version of the regular representation of $\mathfrak{S}_n$.  Using a straightening algorithm on monomials and the Garsia-Stanton basis, Adin, Brenti, and Roichman gave a description of the Frobenius image of $R_n$, graded by partitions, in terms of descents of standard Young tableaux.  Motivated by the Delta Conjecture of Macdonald polynomials, Haglund, Rhoades, and Shimozono gave an extension of the coinvariant algebra $R_{n,k}$ and an extension of the Garsia-Stanton basis.  Chan and Rhoades further extend these results from $\mathfrak{S}_n$ to the complex reflection group $G(r,1,n)$ by defining a $G(r,1,n)$ module $S_{n,k}$ that generalizes the coinvariant algebra for $G(r,1,n)$.  We extend the results of Adin, Brenti, and Roichman to $R_{n,k}$ and $S_{n,k}$ and connect the results for $R_{n,k}$ to skew ribbon tableaux and a crystal structure defined by Benkart et al.

\end{abstract}

\maketitle

\section{Introduction}

The classical coinvariant algebra $R_n$ is constructed as follows: let the symmetric group $\mathfrak{S}_n$ act on the polynomial ring $\mathbb{Q}[x_1,x_2,\ldots, x_n]$ by permutation of the variables $x_1,\ldots, x_n$.  The polynomials that are invariant under this action are called symmetric polynomials, and we let $I_n$ be the ideal generated by symmetric polynomials with vanishing constant term.  Then $R_n$ is defined as the algebra obtained by quotienting $\mathbb{Q}[x_1,x_2,\ldots, x_n]$ by $I_n$, that is 
\begin{equation}R_n:=\frac{\mathbb{Q}[x_1,x_2,\ldots,x_n]}{I_n}.\end{equation}

There are a number of sets of symmetric polynomials in $x_1,\ldots, x_n$ that algebraically generate all symmetric polynomials in the variables $x_1,\ldots, x_n$ with vanishing constant term.  The set that is important for the generalization of $R_n$ that we are considering is the elementary symmetric functions \begin{equation}e_d:=\sum_{1\leq i_1<i_2<\ldots< i_d\leq n} \prod_{j=1}^d x_{i_j},\end{equation}
 for $1\leq d\leq n$.  We then have

 \begin{equation}I_n=\langle e_1,e_2,\ldots, e_n\rangle.\end{equation}

	Since $I_n$ is homogeneous and invariant under the action of $\mathfrak{S}_n$, the coinvariant algebra is a graded $\mathfrak{S}_n$-module.  Since the conjugacy classes of $\mathfrak{S}_n$ are indexed partitions of $n$, the irreducible representations of $\mathfrak{S}_n$ are also indexed by partitions of $n$.  We let $S^\lambda$ denote the irreducible representation corresponding to $\lambda$, and we let $\chi^\lambda_\mu$ be the character of $S^\lambda$ evaluated at an element of type $\mu$.

	The following relies on some definitions that we will cover in Section \ref{Def}.
	Given a representation $V$ of $\mathfrak{S}_n$, a natural question to ask is: "What is the multiplicity of $S^\lambda$ in $V$ for each partition of $n$?".  All of this information can be contained in a single symmetric function called the {\bf Frobenius image} of $V$, which is denoted $Frob(V)$.  The Frobenius image has the following formula
\begin{equation}
	Frob(V)=\sum_{\lambda\vdash n} c_\lambda s_\lambda,
\end{equation}
	where $c_\lambda$ is the multiplicity of $S^\lambda$ in $V$ and $s_\lambda$ is the Schur function associated to $\lambda$.  We will take this formula as a definition.
   In the case of the classical coinvariant algebra this problem was solved by Chevalley \cite{C} who showed that the multiplicity of $S^\lambda$ in $R_n$ is the number of standard Young tableaux of shape $\lambda$, that is that 
\begin{equation}
	Frob(R_n)=\sum_{T\in SYT(n)} s_{sh(T)}
\end{equation}

	If $V$ is a graded representation of $\mathfrak{S}_n$ with degree $d$ component $V_d$, then we can also consider the Frobenius image of $V_d$ for all $d$.  This data can be combined into a single function called the {\bf graded Frobenius image}, which is defined as follows:
\begin{equation}
	grFrob(V;q)=\sum_{d=0}^\infty q^dFrob(V_d).
\end{equation}

	 Lusztig (unpublished) and Stanley \cite{S} showed that for the classical coinvariant algebra the multiplicity of $S^\lambda$ in the degree $d$ component of $R_n$ is the number of standard Young tableaux with major index equal to $d$.  Stated in terms of the graded Frobenius image,
\begin{equation}
grFrob(R_n;q):=\sum_{T\in SYT(n)} q^{maj(T)}s_{shape(T)}.
\end{equation}

	A further refinement of $R_n$ is given as follows: define 

\begin{equation}P_{\trianglelefteq\mu}:=\text{span}\{m\in \mathbb{Q}[x_1,\ldots, x_n]: \lambda(m)\trianglelefteq \mu\},\end{equation}

 and

			\begin{equation}P_{\triangleleft\mu}:=\text{span}\{m\in \mathbb{Q}[x_1,\ldots, x_n]: \lambda(m)\triangleleft \mu\}\end{equation}

		where $m$ are monomials, $\lambda(m)$ is the exponent partition of $m$, and $\triangleleft$ is the dominance order on partitions.  Then let $Q_{\trianglelefteq\mu}$ and $Q_{\triangleleft\mu}$ be the projections of  $P_{\trianglelefteq\mu}$ and $P_{\triangleleft\mu}$ onto $R_n$ respectively.  Next define
		
		\begin{equation}R_{n,\mu}:=Q_{\trianglelefteq\mu}/Q_{\triangleleft\mu}.\end{equation}

  This is a refinement of the grading since the degree $d$ component of $R_n$ is equal to 

		\begin{equation}\bigoplus_{\mu\vdash d} R_{n,\mu}.\end{equation}

	Adin, Brenti, and Roichman \cite{ABR} show that $R_{n,\mu}$ is zero unless $\mu$ is a partition with at most $n-1$ parts such that the differences between consecutive parts are at most 1.  We call such partitions descent partitions.  They also show that in the case that $R_{n,\mu}$ is not zero, the multiplicity of $S^\lambda$ in $R_{n,\mu}$ is given by the number of standard Young tableaux of shape $\lambda$ with descent set equal to the descent set of $\mu$.  Where a descent of a partition $\mu$ is a value $i$ such that $\mu_i>\mu_{i+1}$.  For example if $n=5$, and $\mu=(3,2,2,1)$, then the descents of $\mu$ are $1,3,4$, and the multiplicty of $S^{(2,2,1)}$ is 1 since the only standard Young tableau of shape $(2,2,1)$ with descent set $\{1,3,4\}$ is \begin{equation}\ytableaushort{13,24,5}.\end{equation}

	Motivated by the Delta Conjecture in the theory of Macdonald polynomials, Haglund, Rhoades, and Shimozono \cite{HRS} generalize this entire picture by defining the ideal 

		\begin{equation}I_{n,k}:=\langle x_1^k,x_2^k,\ldots x_n^k, e_{n},e_{n-1},\ldots, e_{n-k+1}\rangle,\end{equation}

for a positive integer $k\leq n$. They then define a generalized coinvariant algebra as 

		\begin{equation}R_{n,k}:=\frac{\mathbb{Q}[x_1,\ldots, x_n]}{I_{n,k}}.\end{equation}

  This is a generalization since in the case $n=k$, we recover the classical coinvariant algebra $R_n$, that is $R_{n,n}=R_n$.  This is connected to the Delta Conjecture because Haglund, Rhoades, and Shimozono show that \begin{equation} (rev_q\circ \omega)\text{grFrob}(R_{n,k};q)\end{equation} is equal to one of the four expressions $\text{Rise}_{n,k}(x;q,0),\text{Rise}_{n,k}(x;0,q),$ $\text{Val}_{n,k}(x;q,0),$ or $\text{Val}_{n,k}(x;0,q)$, where $\text{Rise}_{n,k}$ and $\text{Val}_{n,k}$ are combinatorially defined functions appearing in the Delta Conjecture, and $\omega$ is the standard involution on symmetric functions.

	As in the classical case, $R_{n,k}$ is a graded $\mathfrak{S}_n$-module and we can refine the grading as follows.

	\begin{dfn}
	Let $\mu$ be a partition with at most $n$ parts.  Next define $S_{\trianglelefteq\mu}$ and $S_{\triangleleft\mu}$ to be the projections of $P_{\trianglelefteq\mu}$ and $P_{\triangleleft\mu}$ onto $R_{n,k}$.  We then define 

		\begin{equation}R_{n,k,\mu}:=S_{\trianglelefteq\mu}/S_{\triangleleft\mu}.\end{equation}

	\end{dfn}
	This is a refinement of the grading since the degree $d$ component of $R_{n,k}$ is equal to 

		\begin{equation}\bigoplus_{\mu\vdash d} R_{n,k,\mu}.\end{equation}

Our primary goal is to determine the multiplicties of $S^\lambda$ in $R_{n,k,\mu}$ which we do in the following Theorem, thus extending the results of Adin, Brenti, and Roichman on $R_{n,\mu}$ to $R_{n,k,\mu}$ and refining the results of Haglund, Rhoades and Shimozono.

\begin{thm}
\label{Result1}
	The algebra $R_{n,k,\rho}$ is zero unless $\rho$ fits in an $(n-1)\times k$ rectangle and $\rho_i-\rho_{i+1}\leq 1$ for $i>n-k$.  In the case that $R_{n,k,\rho}$ is not zero, the multiplicity of $S^\lambda$ in $R_{n,k,\rho}$  is given by 

		\begin{equation}\vert \{T\in SYT(\lambda): Des_{n-k+1,n}(\rho)\subseteq Des(T)\subseteq Des(\rho) \}\vert.\end{equation}

\end{thm}  

A key compnent of the methods in \cite{ABR} is the use of a basis for $\mathbb{Q}[x_1,\ldots, x_n]$ that arises from the theory of Cohen-Macaulay rings and the fact that $e_n, e_{n-1},\ldots, e_1$ form a regular sequence.  We are not able to use these methods since the generators of $I_{n,k}$ do not form a regular sequence.

A different direction of generalization comes from considering the coinvariant algebra for general complex reflection groups $G(r,p,n)$, which reduce to $\mathfrak{S}_n$ in the case $r=p=1$.  These algebras are studied by Bagno and Biagioli in \cite{BB}.  Chan and Rhoades \cite{CR} give generalizations of these objects in the case $p=1$ for a parameter $k\leq n$.  The ideal for the case we are interested in is 

		\begin{equation}J_{n,k}:=\langle x_1^{kr},x_2^{kr},\ldots, x_n^{kr}, e_n(\bold{x}_n^r),e_{n-1}(\bold{x}_n^r),\ldots, e_{n-k+1}(\bold{x}_n^r)\rangle,\end{equation}

and the algebra is 

		\begin{equation}S_{n,k}:=\frac{\mathbb{C}[x_1,x_2,\ldots, x_n]}{J_{n,k}}\end{equation}

 where $\bold{x}_n^r,$ denotes the set of variables $\{x_1^r,x_2^r,\ldots, x_n^r\}$.  This is a graded $G(r,1,n)$-module and we can again refine the grading by partitions of size $d$ as follows.

\begin{dfn}

	Let $\mu$ be a partition with at most $n$ parts.  Next define $S_{\trianglelefteq\mu}$ and $S_{\triangleleft\mu}$ to be the projections of $P_{\trianglelefteq\mu}$ and $P_{\triangleleft\mu}$ onto $S_{n,k}$.  We then define 

		\begin{equation}S_{n,k,\mu}:=S_{\trianglelefteq\mu}/S_{\triangleleft\mu}.\end{equation}

\end{dfn}
	This refines the grading since the degree $d$ component of $S_{n,k}$ is equal to 
		
		\begin{equation}\bigoplus_{\mu\vdash d} S_{n,k,\mu}.\end{equation}

The following Theorem gives the multiplicities of irreducible representations appearing in $S_{n,k,\mu}$. 

\begin{thm}
\label{Result2}
	The algebra $S_{n,k,\rho}$ is zero unless $\rho$ fits in an $n\times (kr)$ rectangle, $\rho_i-\rho_{i+1}\leq r$ for $i>n-k$, and $\rho_n<r$.
	In the case that $S_{n,k,\rho}$ is not zero, the multiplicity of $S^{\overline{\lambda}}$ in $S_{n,k,\rho}$  is given by 

		\begin{equation}\vert \{T\in SYT(\overline{\lambda}):
Des^r_{n-k+1,n}(\rho)\subseteq Des(T)\subseteq Des^r(\rho),  c_i(T)\equiv \rho_i (\text{mod }r)\}\vert\end{equation}

\end{thm}

	The paper is organized as follows: Section \ref{Def} will cover background material, Section \ref{Rnk} will prove Theorem \ref{Result1} and show a connection to crystals and skew ribbon tableaux, in Section \ref{Wreath} we will give background for and prove Theorem \ref{Result2}, and in Section \ref{Conclusion} we will cover directions and methods for future work.

\section{Definitions and Background}
\label{Def}

\subsection{Descents and Monomials}

	An important component of the results of \cite{ABR} on $R_n$ is the use of a certain monomial basis for $R_n$.  We will recall this basis and the generalization of this basis given in \cite{HRS} for $R_{n,k}$.  This basis for $R_n$ will be indexed by permutations, and will be defined in terms of the descents of the corresponding permutation.

	Given a permutation $\sigma\in \mathfrak{S}_n$, $i$ is a descent of $\sigma$ if $\sigma(i)>\sigma(i+1)$.  We denote by Des($\sigma$) the set of descents of $\sigma$. We denote by $d_i(\sigma)$, the number of descents of $\sigma$ that are at least as large as $i$, that is 

		\begin{equation}d_i(\sigma):=\vert \{i,i+1,\ldots, n\}\cap Des(\sigma)\vert.\end{equation}

  Finally for two integers $i,j$ such that $1\leq i\leq j\leq n$ we let $Des_{i,j}(\sigma)$ denote the set of descents of $\sigma$ that are between $i$ and $j$ inclusively, that is 

		\begin{equation}Des_{i,j}(\sigma):=Des(\sigma)\cap \{i,i+1,\ldots, j-1, j\}.\end{equation} 

For example if $\sigma=31427865\in \mathfrak{S}_8$, then
\begin{align}
 Des(\sigma)&=\{1,3,6,7\},\\  (d_1(\sigma),\ldots, d_8(\sigma))&=(4,3,3,2,2,2,1,0), \text{and }\\  Des_{2,6}(\sigma)&=\{3,6\}.
\end{align}

Descents are used to define a set of monomials which descend to a basis for $R_n$, see \cite{G} \cite {GS}.

\begin{dfn}

	Given a permutation $\sigma\in\mathfrak{S}_n,$ the {\bf Garsia-Stanton monomial} or simply { \bf descent monomial} associated to $\sigma$ is 

		\begin{equation}gs_\sigma:=\prod_{i=1}^n x_{\sigma(i)}^{d_i(\sigma)}.\end{equation}  These monomials descent to a basis for $R_n$.

\end{dfn}
For example, if $\sigma=31427865\in \mathfrak{S}_8$, then 

		\begin{equation}gs_\sigma=x_3^4x_1^3x_4^3x_2^2x_7^2x_8^2x_6^1\end{equation}

These monomials are generalized by Haglund, Rhodes, and Shimizono in \cite{HRS} to $(n,k)$-descent monomials that are indexed by ordered set partitions of $n$ with $k$ blocks.  Alternatively they can be indexed by pairs $(\pi,I)$ consisting of a permutation $\pi\in\mathfrak{S}_n$ and a sequence $i_1,\ldots, i_{n-k}$ such that 

		\begin{equation}k-des(\pi)>i_1\geq i_2\geq\ldots \geq i_{n-k}\geq 0.\end{equation}

 This is done as follows:

\begin{dfn}

		Given a permutation $\pi\in\mathfrak{S}_n$ and a sequence $I=(i_1,i_2,\ldots, i_{n-k})$ such that 

		\begin{equation}k-des(\pi)>i_1\geq i_2\geq\ldots \geq i_{n-k}\geq 0,\end{equation}

	 the {\bf $(n,k)$-descent monomial} associated to $(\pi,I)$ is

		\begin{equation}gs_{\pi,I}:=gs_\pi x_{\pi(1)}^{i_1}x_{\pi(2)}^{i_2}\ldots x_{\pi(n-k)}^{i_{n-k}}\end{equation}

\end{dfn}

These monomials descend to a basis for $R_{n,k}$.

As an example if $\sigma=31427865\in \mathfrak{S}_8$, $k=6$, and $I=(1,0)$, then\begin{equation}gs_{\sigma, I}=gs_\sigma\cdot x_3^1x_1^0=x_3^5x_1^3x_4^3x_2^2x_7^2x_8^2x_6^1.\end{equation}

\subsection{Permutation and Partitions}

The way that Adin, Brenti, and Roichmann\cite{ABR} make use of the classical descent monomial basis is by using a basis for $\mathbb{Q}[x_1,\ldots, x_n]$ given by Garsia in \cite{G}.  This basis is the set $\{gs_\pi e_\mu\}_{\pi\in\mathfrak{S}_n,\mu\vdash n}$, where 

		\begin{equation}e_\mu=e_{\mu_1}e_{\mu_2}\ldots e_{\mu_{\ell(\mu)}}.\end{equation}

	  In making use of this basis it is necessary to associate certain permutations and partitions to monomials.  Our results also use these, so we recall them here.  

\begin{dfn}

The {\bf index permutation} of a monomial $m=\prod_{i=1}^n x_i^{p_i}$ is the unique permutation $\pi$, such that the following hold:

	\begin{enumerate}

		\item $p_{\pi(i)}\geq p_{\pi(i+1)}$

		\item $p_{\pi(i)}=p_{\pi(i+1)}\implies \pi(i)<\pi(i+1)$

	\end{enumerate}

We denote the index permutation of $m$ as $\pi(m)$.

\end{dfn}

\begin{dfn}

	The {\bf exponent partition} of a monomial $m=\prod_{i=1}^n x_i^{p_i}$ is the partition $(p_{\pi(1)},p_{\pi(2)},\ldots, p_{\pi(n)})$, where $\pi=\pi(m)$.  We denote the exponent partition of $m$ as $\lambda(m)$.

\end{dfn}

We note that if $\lambda$ is the exponent partition of a descent partition, then $\lambda_n=0$ and $\lambda_i-\lambda_{i+1}\leq 1$.  We call a partition that satisfies these conditions a {\bf descent partition}.  If $\lambda$ is the exponent partition of an $(n,k)$-descent monomial, then $\lambda$ has at most $n$ parts, and it has parts of size less than $k$.  We call such partitions {\bf $(n,k)$-partitions}.

\begin{dfn}

	The {\bf complementary partition} of a monomial $m$ is the partition that is conjugate to $(\lambda_i-d_i(\pi))_{i=1}^n$, where $\pi=\pi(m)$ and $\lambda=\lambda(m)$.  We denote the complementary partition of $m$ as $\mu(m)$.

\end{dfn}

To clarify these definitions we present an example.

\begin{ex}

	Let $n=8$, $k=5$, $I=(2,1,1)$ and let 

		\begin{equation}m=x_1^6x_2x_3x_4^2x_6^4x_7x_8^2=x_1^6x_6^4x_4^2x_8^2x_2x_3x_7,\end{equation}

	 then \begin{align} \pi(m)&=16482375,\\  \lambda(m)&=(6,4,2,2,1,1,1,0),\\Des(\pi(m))&=\{2,4,7\},\\ gs_{\pi(m)}&=x_1^3x_6^3x_4^2x_8^2x_2x_3x_7,\\\mu(m)'&=(3,1),\\ \mu(m)&=(2,1,1),\text{and }\\  gs_{\pi(m),I}&=x_1^5x_6^4x_4^3x_8^2x_2x_3x_7.\end{align}

\end{ex}

	The final key component is a partial ordering on monomials of a given degree together with a result on how multiplying monomials by elementary symmetric functions interacts with this partial order.  For a proof of Proposition \ref{Order} we refer the reader to \cite{ABR}.

\begin{dfn}

For $m_1,m_2$ monomials of the same total degree, $m_1\prec m_2$ if one of the following holds:
	\begin{enumerate}
	
		\item $\lambda(m_1)\triangleleft \lambda(m_2)$

		\item $\lambda(m_1)=\lambda(m_2)$ and $inv(\pi(m_1))>inv(\pi(m_2))$
	
	\end{enumerate}
	Where $\triangleleft$ is the strict dominance order on partitions.
\end{dfn}

This partial order is useful because of how it interacts with multiplication of monomials and elementary symmetric functions.  This interaction is encapsulated in the following proposition:

\begin{prop}

\label{Order}

Let $m$ be a monomial equal to $x_1^{p_1}\ldots x_n^{p_n}$, then
among the monomials appearing in  $m\cdot e_{\mu}$,  the monomial 

		\begin{equation}\prod_{i=1}^n x_{\pi(i)}^{p_{(\pi(i))}+u'_i}\end{equation}

 is the maximum with respect to $\prec$, where $\pi$ is the index permutation of $m$.
	\begin{proof}
	We refer the reader to \cite{ABR} for a proof of this theorem.
	\end{proof}

\end{prop}

\subsection{Standard Young Tableaux}

	Our main results come in the form of counting certain standard Young tableaux.

A {\bf Ferrers diagram} is a collection of unit boxes which, since we are using English notation, are justified to the left and up.  The lengths of the rows of a Ferrers diagram form a partition which we call the {\bf shape} of the Ferres diagram.  A {\bf semistandard Young tableau} of size $n$ is a Ferrers diagram containing $n$ boxes where each box is assigned a positive integer such that the intergers increase weakly along rows and strictly down columns.  A {\bf standard Young tableau} is a semistandard Young tableau containing exactly the integers $1,2,\ldots, n$.  We donote the set of standand [semistandard] Young tableaux of size $n$ by $SYT(n) [SSTY(n)]$.  For a partition $\mu$, we let $SYT(\mu) [SSTY(\mu)]$ denote the set of all standard [semistandard] Young tableaux of shape $\mu$.  The weight of a semistandard Young tableau $T$ is the vector $wt(T)$ where the $i$th entry of $wt(T)$ is the number of $i$'s in $T$.  The Schur functions are then defined as 
\begin{equation}
s_\lambda:=\sum_{\lambda\in SSYT(\lambda)} \bold{x}^{wt(\lambda)}
\end{equation}
The Schur functions form a linear basis for symmetric polynomials, and there is a well known involution $\omega$ on the space of symmetric functions that sends $s_\lambda$ to $s_{\lambda'}$ where $\lambda'$ is the partition conjugate to $\lambda$.

An integer $i$ is a {\bf descent} of a standard Young tableaux $T$ if the box containing $i+1$ is strictly below the box containing $i$. We denote by $Des(T)$ the set of all descents of $T$.  Furthermore given two integers $1\leq i\leq j \leq n$ we define $Des_{i,j}$ to be the set of descents of $T$ that are between $i$ and $j$ inclusively, that is\begin{equation}Des_{i,j}(T):=Des(T)\cap\{i,i+1,\ldots j-1,j\}.\end{equation}

\medskip

	As  examples, consider the following Young tableaux:

\medskip

	$T_1=$\ytableaushort{1467,258,3}\,
	$T_2=$\ytableaushort{1347,2568}\,
	$T_3=$\ytableaushort{12478,356}\,

\medskip

	$S_1=$\ytableaushort{1124,233,3}\,
	$S_2=$\ytableaushort{1111,2222}\,
	$S_3=$\ytableaushort{12345,234}\,

	$T_1,T_2,T_3$ are standard Young tableaux, and $S_1,S_2,S_3$ are semistandard Young tableaux.
	The shape of $T_1$ and $S_1$ is $(4,3,1)$, the shape of $T_2$ and $S_2$ is $(4,4)$, and the shape of $(T_3)$ and $S_3$ is $(5,3)$.  The descent sets of the standard Young tableaux are as follows: 

\begin{align}Des(T_1)&=\{1,2,4,7\},\\ Des(T_2)&=\{1,4,7\},\\Des(T_3)&=\{2,4\}.\end{align}

  Next, $Des_{5,7}(T_1)=Des_{5,7}(T_2)=\{7\}$, and $Des_{5,7}(T_3)=\emptyset$
The weight of the semistandard Young tableaux are as follows:

\begin{align}wt(S_1)&=(2,2,4,1,0,0,\ldots),\\ wt(S_2)&=(4,4,0,0,0,0\ldots),\\wt(S_3)&=(1,2,2,2,1,0,\ldots).\end{align}

A {\bf skew Young tableau} is a Young tableau that has had a Young tableau removed from its upper left corner.  The definitions of both semistandard Young tableaux and Schur function extend to {\bf semistandard skew Young tableaux} and {\bf skew Schur functions}.  A connected skew Young tableau that does not contain any $2\times 2$ boxes is called a skew ribbon tableau.  These two conditions make it so that the shape of a skew ribbon tableau is uniquely determined by the lengths of its rows, so that we can specify a skew-ribbon tableau shape by a sequence of postive integers.  For example if we specify that a skew ribbon tableau has rows of lengths $(4,2,1,3)$, then the following are two examples of semistandard skew Young tableaux with the only possible shape:

	\ytableaushort{\none\none\none1123,\none\none12,\none\none3,124}\,\,\,
	\ytableaushort{\none\none\none1111,\none\none12,\none\none2,113}\,.

\section{Descent representations of $R_{n,k}$}
\label{Rnk}
In the case of the classical coinvariant algebra, Adin, Brenti, and Roichman determine the isomorphism type of $R_{n,\rho}$ by comparing the graded traces of the actions of $\mathfrak{S}_n$ on $\mathbb{Q}[x_1,\ldots, x_n]$  and on $R_{n}$.  We will follow a similar path, but instead of considering the action of $\mathfrak{S}_n$ on $\mathbb{Q}[x_1,\ldots, x_n]$, we will consider its action on the space 

		\begin{equation}P_{n,k}:=\text{span}_{\mathbb{Q}}\{ x_1^{p_1}x_2^{p_2}\ldots x_n^{p_n}:p_1,p_2,\ldots, p_n<k\},\end{equation}

that is the space of rational polynomials in the variables $x_1,\ldots, x_n$ where the powers of each $x_i$ are less than $k$.

We begin by giving a straightening Lemma that is a similar to a Lemma of Adin, Brenit, and Roichmann \cite{ABR}.  Our Lemma differs from theirs in that we are considering monomials in $P_{n,k}$ instead of $\mathbb{Q}[x_1,\ldots, x_n]$, we use $(n,k)$-descent monomials instead of the classical descent monomials, and we consider elementary symmetric functions corresponding to partitions with parts of size at least $n-k+1$ instead of all elementary symmetric functions.
\begin{lem}

\label{OrderLem}
If  $m=\prod_{i=1}^n x_i^{p_i}$ is a monomial in $P_{n,k}$ (that is $p_i<k$ for all $i$), then 

		\begin{equation}m=gs_{\pi,I}e_\nu+\sum.\end{equation}

  Where $\pi=\pi(m)$; $\sum$ is a sum of monomials $m'\prec m$; I is the length $n-k$ sequence defined by $i_\ell=\mu'_{\ell}-\mu'_{n-k+1}$, where $\mu$ is the complementary partition of $m$; and $\nu$ is the partition specified by:

	\begin{enumerate}

		\item $\nu'_\ell=\mu'_\ell$ for $\ell> n-k$

		\item $\nu'_\ell=\mu'_{n-k+1}$ for $\ell\leq n-k$
	
	\end{enumerate}

Furthermore $\nu$ consists of parts of size at least $n-k+1$.		

	\begin{proof}
	In order to show that $gs_{\pi,I}$ is well defined we need to check that $k-des(\pi)>i_1\geq i_2\geq\ldots\geq i_k\geq 0$.  By definition $i_1=\mu'_1-\mu'_{n-k+1}\leq \mu'_1=p_{\pi(1)}-d_1(\pi)$ and then by assumption $p_{\pi(1)}<k$, and $d_1(\pi)=des(\pi)$, thus 

		\begin{equation}i_1\leq p_{\pi(1)}-d_1(\pi)<k-des(\pi).\end{equation}

	We also note that I is a non-negative weakly-decreasing sequence since it consists of the parts of a partition minus a constant that is at most as large as the smallest part of the partition.  Thus $I$ satifies the condition so $gs_{\pi,I}$ is well defined.
	
	Next we show that $gs_{\pi,I}$ and $m$ have the same index permutation, that is that 

		\begin{equation}\pi(gs_{\pi,I})=\pi(m)=\pi.\end{equation}

  To show this, we need to consider the sequence of the exponents of $x_{\pi(\ell)}$ in $gs_{\pi,I}$.  This sequence is the sum of the sequences $d_\ell(\pi)$ and $i_\ell$ (where we take $i_\ell=0$ for $\ell>n-k$).  Since these are both weakly-decreasing sequences, their sum is also weakly-decreasing.  Furthermore if $d_\ell(\pi)+i_\ell=d_{\ell+1}(\pi)+i_{\ell+1}$, then $d_\ell(\pi)=d_{\ell+1}(\pi)$, which by the definition of $d_\ell(\pi)$ implies that $\ell$ is not a descent of $\pi$, that is that $\pi(\ell)<\pi(\ell+1)$, thus $\pi$ satisfies the two conditions of being the index permutation, and thus by uniqueness it is the index permutation.

	Now by Proposition \ref{Order}, the maximum monomial in $gs_{\pi,I}e_\nu$ will have the form $\prod_{\ell=1}^n x_{\pi(\ell)}^{q_{\ell}}$ where $q_\ell$ is given by:

	\begin{enumerate}

		\item $q_\ell=d_\ell(\pi)+i_\ell+\nu'_\ell$ for $\ell\leq n-k$

		\item $q_\ell=d_\ell(\pi)+\nu'_\ell$ for $\ell> n-k$

	\end{enumerate}

By substitution, first using the definitions of $i_\ell$ and $\nu_\ell$ and then the definition of the complementary partition, we get 

		\begin{equation}q_\ell=d_\ell(\pi)+\mu'_\ell-\mu'_{n-k+1}+\mu'_{n-k+1}=d_\ell(\pi)+\mu'_\ell=p_{\pi(\ell)}\end{equation}

 for $\ell\leq n-k$, and 

		\begin{equation}q_\ell=d_\ell(\pi)-\mu'_\ell=p_{\pi(\ell)}\end{equation}

for $\ell>n-k$.

Finally, $\nu$ has parts of size at least $n-k+1$ because by definition, the first $n-k+1$ parts of $\nu'$ are all the same size.

	\end{proof}
\end{lem}

This Lemma gives rise to a basis for $P_{n,k}$ which will be key to relating how $\mathfrak{S}_n$ acts on $P_{n,k}$ to how it acts on $R_{n,k}$.

\begin{prop}
\label{BasisProp}

	The set $D_{n,k}$ consisting of products $gs_{\pi,I}e_\nu$ where $\nu$ is a partition with parts of size at least $n-k+1$ and $(\lambda(gs_{\pi,I})+\nu')_1<k$ form a basis for $P_{n,k}$.
	\begin{proof}

The condition that $(\lambda(gs_{\pi,I})+\nu')_1<k$ along with Lemma \ref{OrderLem} guarantees that the maximum monomial in each element of $D_{n,k}$ is contained in $P_{n,k}$.  Then since the partial order $\prec$ refines dominace order, all other monomials appearing in elements in $D_{n,k}$ are also contained in $P_{n,k}$.  Therefore $D_{n,k}$ is contained in $P_{n,k}$.  

Iteratively applying Lemma \ref{OrderLem} lets us express any monomial in $P_{n,k}$ as a linear combination of elements in $D_{n,k}$, thus $D_{n,k}$ spans $P_{n,k}$.    To show that this expansion is unique (up to rearrangement) it is sufficient to show that if the maximal monomials in $gs_{\pi,I}e_{\nu}$ and $gs_{\phi,J}e_{\rho}$ are the same, then $\pi=\phi$, $I=J$, and $\nu=\rho$.  To see this, we note that as a corollary of the proof of Lemma \ref{OrderLem}, the index permutations of the maximal monomials are the same, and they are $\pi$ and $\phi$ respectively, and thus $\pi=\phi$.  Then, by Proposition \ref{Order}, the power of $x_{\pi(\ell)}$ in each of these maximum monomials will be $d_\ell(\pi)+i_\ell+\nu'_\ell$ and $d_\ell(\pi)+j_\ell+\rho'_\ell$.  This immediately gives that $\nu'_\ell=\rho'_\ell$ for $\ell>n-k$ since $i_\ell=j_\ell=0$ for $\ell>n-k$.  Then since the first $n-k+1$ parts of $\nu'$  are all equal and the first $n-k+1$ parts of $\rho'$ are equal and since $\nu'_{n-k+1}=\rho'_{n-k+1}$, we have that $\nu'=\rho'$ which implies $\nu=\rho$.  This then implies that $i_\ell=j_\ell$ for all $\ell$, and therefore this expansion is unique.  Therefore $D_{n,k}$ is linearly independent and is a basis.

	\end{proof}

\end{prop}

\begin{prop}
\label{DomProp}
	Let $p$ be the map projecting from $\mathbb{Q}[x_1,\ldots, x_n]$ to $R_{n,k}$ and let $m$ be a monomial in $P_{n,k}$. Then 

		\begin{equation}p(m)=\sum_{\pi,I} \alpha_{\pi,I} gs_{\pi,I}\end{equation}

	 where $\alpha_{\pi,I}$ are some constants, and the sum is over pairs $(\pi,I)$ such that $\lambda(gs_{\pi,I})\trianglelefteq \lambda(m)$.

	\begin{proof}
		Since $D_{n,k}$ is a basis, we can express $m=\sum_{\pi,I,\nu} \alpha_{\pi,I,\nu} gs_{\pi,I} e_\nu$ for some constants $\alpha_{\pi,I\nu}$.  By Lemma \ref{OrderLem}, $\alpha_{\pi,I,\nu}$ is zero if the leading monomial of $gs_{\pi,I} e_\nu$ is not weakly smaller than $m$ under the partial order on monomials.  But since the partial order on monomials refines the dominance order on exponent partitions, for each non-zero term the exponent partition of the leading monomial will be dominated by $\lambda(m)$, that is that 

		\begin{equation}(\lambda(gs_{\pi,I})+\nu')\trianglelefteq \lambda(m).\end{equation}

	  Then when we project down to $R_{n,k}$, each term with $\nu\neq \emptyset$ will vanish since $e_\nu$ is in $I_{n,k}$, so that, 

		\begin{equation}p(m)=\sum_{\pi,I}\alpha_{\pi,I,\emptyset} gs_{\pi,I}\end{equation}

 where the sum is over $(\pi,I)$ such that $\lambda(gs_{\pi,I})\trianglelefteq \lambda(m)$.
	\end{proof}

\end{prop}

This Proposition gives the following Corollary:

\begin{coro}
\label{zero}

	$R_{n,k,\rho}$ is zero unless $\rho$ is the exponent partition of an $(n,k)$-descent monomial, which occurs precisely when $\rho$ is an $(n,k)$ partition such that the last $k$ parts form a descent partition.

\end{coro}

This basis allows us to express the action of $\tau\in\mathfrak{S}_n$ on $ P_{n,k}$ in terms of its action on $R_{n,k}$ with the basis of $(n,k)$-Garsia-Stanton monomials.  That is if\begin{equation}\tau(gs_{\pi,I})=\sum_{\phi,J} \alpha_{\phi,J} gs_{\phi,J},\end{equation}for some constants $\alpha_{\phi,J}$, then\begin{equation}\tau(gs_{\pi,I} e_\nu)=\sum_{\phi,J} \alpha_{\phi,J} gs_{\phi,J} e_\nu.\end{equation}

This equality holds since $e_\nu$ is invariant under the action of $\mathfrak{S}_n$.  The important part is that for each $(\phi,J)$ such that $\alpha_{\phi,J}$ is non-zero, $gs_{\phi,J}e_\nu$ is in $D_{n,k}$, that is we do not violate the condition that $(\lambda(gs_{\phi,J})+\nu')_1<k$.  By Proposition \ref{DomProp} we have that $\lambda(gs_{\phi,J})\triangleleft \lambda(gs_{\sigma,I})$ which implies that $\lambda(gs_{\phi,J})+\nu'\triangleleft \lambda(gs_{\sigma,I})+\nu'$, and thus $(\lambda(gs_{\phi,J})+\nu')_1\leq (\lambda(gs_{\sigma,I})+\nu')_1<k$, which implies that we satisfy the condition to be in our basis.

We now move to the Lemmas that will allow us to prove our main result.
	
\begin{lem}

\label{DecomLem}
 Given an $(n,k)$-partition $\mu$ and an $(n,k)$-descent partition $\nu$ there exists a $(n,k)$-partition $\rho$ such that $\mu=\nu+\rho$ if and only if $Des(\nu)\subseteq Des(\mu)$.  If it exists, $\rho$ is unique.

	\begin{proof}
		There is only one possible value for each part of $\rho$ which is $\rho_i=\mu_i-\nu_i$, the only thing to check is whether this gives a partition, specifically we need to check whether $\rho_i-\rho_{i+1}=(\mu_i-\mu_{i+1})-(\nu_i-\nu_{i+1})\geq 0$.  Since $\nu$ is a descent partition, $(\nu_i-\nu_{i+1})$ is 1 if $i$ is a descent of $\nu$ and 0 if it is not.  Similarly, $(\mu_i-\mu_{i+1})$ is at least 1 if $i$ is a descent of $\mu$ and 0 otherwise.  Thus in order for $(\mu_i-\mu_{i+1})-(\nu_i-\nu_{i+1})$ to be non-negative, it is necessary and sufficient that if $i$ is a descent of $\nu$, then $i$ is also a descent of $\mu$.  That is, $\rho$ will be a partition if and only if $Des(\nu)\subset Des(\mu)$.
	\end{proof}
\end{lem}

\begin{ex}
As an example of the Lemma \ref{DecomLem}, let $n=8,k=6$ then let \begin{align}\mu&=(5,5,3,3,1,1,1,0),\\ \nu_1&=(2,2,1,1,0,0,0,0),\\\nu_2&=(3,3,2,2,1,1,0,0).\end{align}

Then \begin{align}Des(\mu)&=\{2,4,7\},\\Des(\nu_1)&=\{2,4\},\\Des(\nu_2)&=\{2,4,6,8\}.\end{align}

We then have that $Des(\nu_1)\subseteq Des(\mu)$, and that $\mu-\nu_1=(3,3,2,2,1,1,1,0)$ is a partition.  On the other hand, $Des(\nu_2)\not\subseteq Des(\mu)$, and $\mu-\nu_2=(2,2,1,1,0,0,1,0)$ is not a partition.

\end{ex}

\begin{lem}

\label{Decom2Lem}

Given an $(n,k)$-partition $\mu$ and a set $S\subseteq Des_{n-k+1,n}(\mu)$, there is a unique pair $(\nu,\rho)$ such that $\mu=\nu+\rho$ and $\nu$ is the exponent partition of an $(n,k)$-descent monomial with $Des_{n-k+1,n}(\nu)=S$, and $\rho$ is an $(n,k)$-partition with $\rho_1=\rho_2=\ldots =\rho_{n-k+1}$.

	\begin{proof}

		The last $k$ values of the exponent partition of a descent monomial form a descent partition, so applying Lemma \ref{DecomLem} to the partition determined by $S$ determines the last $k$ values of $\rho$.   Then since we need that the first $n-k+1$ values of $\rho$ are the same, this determines what $\rho$ must be, and by subtraction what $\nu$ must be.  We just need to check that $\nu$ is actually a partition, that is that $\nu_i-\nu_{i+1}\geq 0$ for $1\leq i\leq n-k$.  This is true since $\nu_i-\nu_{i+1}=\mu_{i}-\mu_{i+1}\geq 0$ because $\rho_i=\rho_{i+1}$ for $i\leq n-k$.

	\end{proof}

\end{lem}

We give an example of how Lemma \ref{Decom2Lem} works.

\begin{ex}

	Let $n=8,k=6$, and let\begin{equation}\mu=(5,5,3,3,1,1,1,0),\end{equation}and let $S=\{4\}$, then\begin{equation}\nu=(4,4,2,2,1,1,1,0),\end{equation}and\begin{equation}\rho=(1,1,1,1,0,0,0,0)\end{equation}

\end{ex}

We now give a proof of Theorem \ref{Result1}
	\begin{proof}(Thm. \ref{Result1})
		The determination of when $R_{n,k,\rho}$ is zero is from Corollary \ref{zero}.
		
		Next we define an inner product on polynomials by $\langle m_1, m_2 \rangle=\delta_{m_1m_2}$ for two monomials $m_1$, $m_2$, and then extending bilinearly.
		 We then consider the graded trace of the action of $\tau\in \mathfrak{S}_n$ on $P_{n,k}$ defined for the monomial basis by\begin{equation}Tr_{P_{n,k}}(\tau):=\sum_m \langle \tau(m),m\rangle \cdot \bar{q}^{\lambda(m)}\end{equation}where $\overline{q}^{\lambda}=\prod_{i=1}^n q_i^{\lambda_i}$ for any partition $\lambda$. Adin, Brenti, Roichman show that
\begin{equation}Tr_{\mathbb{Q}[x_1,\ldots, x_n]}(\tau)=\sum_{\lambda\vdash n} \chi^\lambda_\mu\frac{\sum_{T\in SYT(\lambda)} \prod_{i=1}^n q_i^{d_i(T)}}{\prod_{i=1}^n (1-q_1q_2\ldots q_i)}\end{equation}(where $\mu$ is the cycle type of $\tau$).
From this we can recover $Tr_{P_{n,k}}(\tau)$ by restricting to powers of $q_1$ that are at most $k-1$.  Doing this gives\begin{equation}\sum_{\lambda\vdash n}\chi_\mu^\lambda\sum_{T\in SYT(\lambda),\nu}\bar{q}^{\lambda_{Des(T)}}\bar{q}^\nu.\end{equation} Where the $\nu$'s are partitions such that $(\lambda_{Des(T)})_1+\nu_1<k$, and $\lambda_{Des(T)}$ is the descent partition with descent set $T$.

		Alternatively, we can calculate $Tr_{P_{n,k}}(\tau)$ by using the basis from Proposition \ref{BasisProp}, this gives
\begin{align}Tr_{P_{n,k}}(\tau)&=\sum_{\sigma,I,\nu} \langle \tau(gs_{\sigma,I} e_\nu),gs_{\sigma,I} e_\nu\rangle \bar{q}^{\lambda(gs_{\sigma,I})}\bar{q}^{\nu'}\\
&=\sum_{\sigma,I,\nu} \langle \tau(gs_{\sigma,I}),gs_{\sigma,I} \rangle \bar{q}^{\lambda(gs_{\sigma,I})}\bar{q}^{\nu'}\\
&=\sum_{\lambda,\nu} Tr_{R_{n,k}}(\tau; \bar{q}^\lambda)\bar{q}^{\lambda}\bar{q}^{\nu'}\end{align}
where the $\nu$'s are partitions with parts of size at least $n-k+1$ such that $(\lambda(gs_{\sigma,I}))_1+(\nu')_1<k$, and $Tr_{R_{n,k}}(\tau;\bar{q}^\lambda)$ is the coefficient of $\bar{q}^\lambda$ in the graded trace of the action of $\tau$ on $R_{n,k}$.

		We now consider the coefficient of $\bar{q}^\rho$ for some partition $\rho$.  Using the first calculation and Lemma \ref{DecomLem}, the inner sum can be reduced to $T$ such that $Des(T)\subseteq Des(\rho)$, so that we get\begin{equation}\sum_{\lambda\vdash n} \chi^\lambda_\mu\vert \{T\in SYT(\lambda), Des(T)\subseteq Des(\rho)\}\vert.\end{equation}

Looking at the second calculation and using Lemma \ref{Decom2Lem} gives\begin{equation}\sum_{S\subseteq Des_{n-k+1,n}(\rho)} Tr_{R_{n,k}}(\tau;\bar{q}^{\lambda_S}),\end{equation}where $\lambda_S$ is the exponent partition of some $(n,k)$-descent monomial $gs_{\sigma, I}$ with $Des_{n-k+1,n}(\lambda(gs_{\sigma, I}))=S$.  Together this gives that\begin{equation}\sum_{\lambda\vdash n} \chi^\lambda_\mu\vert \{T\in SYT(\lambda): Des(T)\subseteq Des(\rho)\}\vert=\sum_{S\subseteq Des_{n-k+1,n}(\rho)} Tr_{R_{n,k}}(\tau;\bar{q}^{\lambda_S}).\end{equation}

	We want to further refine this result by showing that\begin{equation}\sum_{\lambda\vdash n} \chi^\lambda_\mu\vert \{T\in SYT(\lambda): S\subseteq Des(T)\subseteq Des(\rho)\}\vert= Tr_{R_{n,k}}(\tau;\bar{q}^{\lambda_S})\end{equation}for any specific $S'$.  We do this by induction on $\vert \lambda_{S'}\vert$.  The base case is $\lambda_{S'}=\emptyset$ is trivial.  If we take $\rho=\lambda_{S'}$, then $\lambda_{S'}$ will appear in the sum since we can take the $\nu$ from Lemma \ref{Decom2Lem} to be 0, and all other $\lambda_S$'s will be smaller since the corresponding $\nu$'s will be non-empty.  Thus by the inductive hypothesis,\begin{equation}\sum_{\lambda\vdash n} \chi^\lambda_\mu\vert \{T\in SYT(\lambda): S'\not\subseteq Des(T)\subseteq Des(\rho) \}\vert= \sum_{S\subsetneq S'}Tr_{R_{n,k}}(\tau;\bar{q}^{\lambda_S}),\end{equation}subtracting this from our result gives the desired refinement.  This then proves the Theorem since the exponent partition of any $(n,k)$-descent monomial $gs_{\sigma, I}$ appears when we take  $\rho=\lambda(gs_{\sigma,I})$.

	\end{proof}

\begin{ex}
Let
$n=8$, $k=6$, $\rho=(5,3,2,2,1,1,1),$ $\lambda=(4,3,1),$
then
$Des_{1,2}(\rho)=\{1,2\},$ and  $Des_{3,8}=\{4,7\}.$

The standard Young tableaux $T$ of shape $\lambda$ with $\{4,7\}\subseteq Des(T)\subseteq \{1,2,4,7\}$ are as follows:

\medskip
	
	\ytableaushort{1467,258,3}*[*(blue!20)]{2,1}*[*(white)]{3}*[*(blue!20)]{4}\,
	\ytableaushort{1347,268,5}*[*(blue!20)]{1}*[*(white)]{2}*[*(blue!20)]{4}\,
	\ytableaushort{1347,256,8}*[*(blue!20)]{1}*[*(white)]{2}*[*(blue!20)]{4}\,

	\medskip

	\ytableaushort{1247,356,8}*[*(white)]{1}*[*(blue!20)]{4}\,
	\ytableaushort{1247,368,5}*[*(white)]{1}*[*(blue!20)]{4}\,
	\ytableaushort{1267,348,5}*[*(white!20)]{1,1}*[*(blue!20)]{2}*[*(white!20)]{3}*[*(blue!20)]{4,2}\,
	\ytableaushort{1234,567,8}*[*(white)]{3,2}*[*(blue!20)]{4,3}

Therefore by Theorem \ref{Result1}, the coefficient of $S^\lambda$ in $R_{n,k,\rho}$ is 7.

\end{ex}

Theorem \ref{Result1} is related to the a crystal structure that defined by Benkart, Colmenarejo, Harris, Orellana, Panova, Schilling, and Yip \cite{Sh}.  Like $R_{n,k}$, the crystal structure that they define is motivated by the Delta Conjecture, and its graded character is equal to \begin{equation}(rev_q\circ \omega) \text{grFrob}(R_{n,k};q),\end{equation} which, as we mentioned before, is equal to a special case of the combinatorial side of the Delta Conjecture.  This crystal is built up from crystal structures on ordered multiset partitions in minimaj ordering with specificed descents sets, and the characters of these smaller crystals is given in terms of skew ribbon tableaux.  Since there is an algebra corresponding to the entire crystal structure, it is natural to wonder if there algebras that correspond to these smaller crystals.  The algebras $R_{n,k,\rho}$ are these algebras.

In order to see this connection, we need to rewrite the Frobenius image of $R_{n,k,\rho}$ that we get from Theorem \ref{Result1} to get an expression in terms of skew-ribbon tableuax.  Using the combinatorial definition of $s_\lambda$ and Theorem \ref{Result1}, we can write Frobenius image of $R_{n,k,\rho}$ as \begin{equation}\text{Frob}(R_{n,k,\rho})=\sum_{(P,Q)} \bold{x}^{wt(P)}\end{equation} where the sum is over pairs $(P,Q)$ with the following conditions:
\begin{itemize}

\item $P$ is a semistandard Young tableau of size $n$
\item $Q$ is a standard Young tableau of size $n$
\item $sh(P)=sh(Q)$
\item $Des_{n-k+1,n}(\rho)\subseteq Des(Q)\subseteq Des(\rho)$.\end{itemize}  

The Robinson-Schensted-Knuth(RSK) correspondence gives a weight-preserving bijection between pairs $(P,Q)$ with the above conditions and words $w$ of length $n$ in the alphabet of positive integers with $Des_{n-k+1,n}(\rho)\subseteq Des(w)\subseteq Des(\rho)$.  Therefore if we apply the reverse RSK correspondence to the Frobenius image it can be rewritten as
 \begin{equation}\text{Frob}(R_{n,k,\rho})=\sum_{w} \bold{x}^{wt(w)}\end{equation} where the sum is over words of length $n$ with $Des_{n-k+1,n}(\rho)\subseteq Des(w)\subseteq Des(\rho)$.

Next let $d_i$ be the difference between the $i$th and $(i-1)$th descents of $\rho$, taking $d_1$ to be the first descent.  Then let $p$ be the index of the largest descent smaller than $n-k+1$.  We this notation any word $w$ as above can be split into subwords $w_1,w_2,\ldots, w_p,$ and $v$ such that $w=w_1w_2\ldots w_pv$ where each $w_i$ has length $d_i$ and has no descents, and $v$ has descents at $d_{p+1}, d_{p+1}+d_{p+2},\ldots, d_{p+1}+d_{p+2}+\ldots d_{des(\rho)}$.  Any such collection of subwords gives an acceptable word $w$, thus $\text{Frob}(R_{n,k,\rho}$ can be written as product of terms of the form $\sum_{w_i} \bold{x}^{wt(w_i)}$ and $\sum_{v} \bold{x}^{wt(v)}$, where the sums are over words with the corresponding restrictions.  These terms can be simplfied as follows.  The term $\sum_{w_i} \bold{x}^{wt(w_i)}$ is equal to $h_{d_i}$ and $\sum_{v} \bold{x}^{wt(v)}$ is equal to $s_{\gamma}$ where $\gamma$ is the skew ribbon shape with rows of lengths $(n-(d_1+d_2+\ldots+d_{des(\rho)}),d_{des(\rho)},d_{des(\rho)-1},\ldots, d_{p+1})$.  This last part is because there is a bijection between fillings of $\gamma$ and words with the conditions of $v$ given by reading the fillings of $\gamma$ row by row from bottom to top reading each row from left to right.  Combining these together gives that the Frobenius image of $R_{n,k,\rho}$ is equal to \begin{equation}\text{Frob}(R_{n,k,\rho})=s_\gamma \prod_{i=1}^p h_{d_i}\end{equation}

To clarify the above we will give an example.  Let $n=11, k=8$, let $\rho=(7,7,5,3,3,3,3,2,1,1)$.  Then $Des_{n-k+1,n}(\rho)=\{7,8,10\}$ and $Des(\rho)=\{2,3,7,8,10\}$, and the values of $d_i$ written in a list are $2,1,4,1,2$ and $p=2$.  The skew ribbon tableau that will appear will thus have row lengths $(1,2,1,4).$  A pair $(P,Q)$ with the above conditions would then be 

\medskip

	P=\ytableaushort{11235,2355,34}\,
	Q=\ytableaushort{12567,348{10},9{11}}*[*(white)]{1}*[*(blue!20)]{2,}*[*(white)]{4,2}*[*(blue!20)]{5,4}\,

\medskip

Applying the reverse RSK correspondence to this pair gives the word \begin{equation}w=34125553132.\end{equation}  This is then broken up into the words $w_1=34$, $w_2=1$, and $v=25553132$ these are then put into semistandard Young (skew) tableaux as follows 

	\begin{center}\ytableaushort{34}\, , \ytableaushort{1}\, , \ytableaushort{\none\none\none\none 2,\none\none\none 13,\none\none\none 3,2555}
	\end{center}

If we apply $\omega$ to this product we get \begin{equation}\omega(\text{Frob}(R_{n,k,\rho}))=s_{\gamma'} \prod_{i=1}^p e_{d_i}.\end{equation}  This expression (for the appropriately chosen values) is the character of the crystals that Benkart, Colmenarejo, Harris, Orellana, Panova, Schilling, and Yip \cite{Sh} use to build up their main crystal strucuture.    Therefore the algebras $R_{n,k,\rho}$ fill in a piece that was missing on the algebraic side of things.

Using Theorem \ref{Result1} we can also recover a result of Haglund, Rhoades, and Shimozono \cite{HRS}.

We will use the $q$-binomial coefficient which has the following formulation.

\begin{equation}[n]_q:=1+q+\ldots q^{n-1}\quad [n]!q:=[n]_q[n-1]_q\ldots [1]_q\quad \genfrac{[}{]}{0pt}{}{n}{k}_q:=\frac{[n]_q}{[k]_q[n-k]_q}\end{equation}

Additionally we will use the well known result that the coefficient of $q^d$ in $\genfrac{[}{]}{0pt}{}{n+m}{m}_q$ is the number of partitions of size $d$ that fit in an $n\times m$ box.

\begin{coro}
	Let $f_\lambda(q)$ be the generating function for the multiplicities of $S^{\lambda}$ in the degree $d$ compent of $R_{n,k}$.  Then\begin{equation}f_\lambda(q)=\sum_{T\in SYT(\lambda)} q^{maj(T)}\genfrac{[}{]}{0pt}{}{n-des(T)-1}{n-k}_q.\end{equation} Where the major index $maj(T)$ is the sum of the descents of $T$.
	\begin{proof}

		By Theorem \ref{Result1}, each standard Young tableau of shape $\lambda$ contributes to $f_\lambda(q)$ once for each partition $\rho$ such that $\rho$ is the exponent partition of an $(n,k)$-descent monomial and $Des_{n-k+1,n}(\rho)\subseteq Des(T)\subseteq Des(\rho)$.  All such $\rho$ come from $(n,k)$-descent monomials $gs_{\pi,I}$ where $\pi$ is a permutation with $Des(\pi)=Des(T)$ and $I$ is a sequence such that $k-des(T)>i_1\geq i_2\geq\ldots\geq i_{n-k}\geq 0$.  This choice of $I$ is the same as choosing a partition that fits in an $(n-k)\times (k-1-des(T))$ box.  The generating function for the number of partitions of size $d$ that fit in an $(n-k)\times (k-1-des(T))$ box is $\genfrac{[}{]}{0pt}{}{(n-k)+(k-des(T)-1)}{n-k}_q=\genfrac{[}{]}{0pt}{}{n-des(T)-1}{n-k}_q$.  The factor of $gs_\pi$ in the $(n,k)$ descent monomial then has degree $maj(T)$, so that each standard Young tableau $T$ of shape $\lambda$ will contribute $q^{maj(T)}\genfrac{[}{]}{0pt}{}{n-des(T)-1}{n-k}_q$  to $f_\lambda(q)$.  This completes the proof.

	\end{proof}
\end{coro}

The proof of this result in \cite{HRS} is fairly involved using a tricky recursive argument involving an auxillary family of algebras.  Our  method gives a simpler proof for the result.

\section{Wreath Products}
\label{Wreath}

This picture can be extended by looking at reflection groups other than $\mathfrak{S}_n$.  Specifically we will look at the complex reflection group $G(r,1,n)$ which is equal to the wreath product of $\mathbb{Z}_r$ and $\mathfrak{S}_n$.  This group acts on $\mathbb{C}[x_1,\ldots, x_n]$ by $\mathfrak{S}_n$ permuting the variables and by the $i$th copy of $\mathbb{Z}_r$ sending $x_i$ to $\xi x_i$ where $\xi$ is a primitive $r$th root of unity.  Alternatively, we can view this group as the set of $n\times n$ matrices with exactly 1 non-zero entry in each row and column where the non-zero entries are $r$th roots of unity.  The action of $G(r,1,n)$ on $\mathbb{C}[x_1,\ldots, x_n]$ in this case is matrix multiplication.  A third way of thinking of this group is as permutation of $n$ in which each number is assigned one out of $r$ colors. 

Throughout this section many of the objects we consider will depend on the positive integer $r$, but since we only ever consider a fixed $r$ we will frequently supresses the $r$ in our notation in order to avoid cumbersome notation.  To begin we will write $G_n$ for the group $G(r,1,n)$.

As in the case of $\mathfrak{S}_n$ there is a coinvariant algebra $S_n$ associated to this action of $G_n$ that is defined as\begin{equation}S_n:=\frac{\mathbb{C}[x_1,\ldots, x_n]}{J_n},\end{equation}where $J_n$ is the ideal generated by all polynomials invariant under the action of $G_n$ with zero constant term.  Any polynomial that is invariant under the action of $G_n$ must be a symmetric polynomials in the variables $x_1^r,x_2^r,\ldots, x_n^r$.  We denote this set of variables as $\bold{x}_n^r$.,  Then $J_n=\langle e_n(\bold{x}_n^r),\ldots e_1(\bold{x}_n^r)\rangle$.

Our goal is to give the multiplicities of all irreducible representations of $G_n$ in $S_{n,k,\rho}$.  In order to do this, we will review some of the representation theory of $G_n$.

\subsection{Backgrounds and Definitions}

The elements of $G_n$ can be viewed as r-colored permutations of length n which are defined as follows:

	\begin{dfn}

		An {\bf r-colored permutatio} of length $n$ is a permutation $\pi=\pi_1\ldots \pi_n$ where each value $\pi_i$ has been assigned a value $c_i$ from the set $\{0,1,\ldots, r-1\}$.  We can write this in the form $\pi_1^{c_1}\pi_2^{c_2}\ldots \pi_n^{c_n}$.
		
	\end{dfn}

	For example $3^05^22^11^24^0$ is a $3$-colored permutation of length $5$.

	As before we define a statistic on r-colored permutation called descents.

	\begin{dfn}
	An index $i$ is a {\bf descent} of an r-colored permutation $g=\pi_1^{c_1}\pi_2^{c_2}\ldots \pi_n^{c_n}$ if one of the following conditions hold: 
\begin{enumerate}
\item $c_i<c_{i+1}$
\item $c_i=c_{i+1}$ and $\pi_i>\pi_{i+1}$.
\end{enumerate}

  We denote the set of descents of an r-colored permutation by $Des(g)$.  Furthermore we will denote $\vert Des(g)\vert$ as $des(g)$, and we will write $d_i(g)$ to be the number of descents of $g$ that are at least as large as $i$, that is\begin{equation}d_i(g):=\vert Des(g)\cap \{i,i+1,\ldots, n\}\vert.\end{equation}

	\end{dfn}

	For example if $g=3^05^24^26^01^12^1$, then $Des(g)=\{1,2,4\}$ since $1$ and $4$ satisfy condition (1) and $2$ satifies condition $(2)$.

	Using these $d_i$ values we follow Bango and Biagoli\cite{BB} in defining flag descent values as\begin{equation}f_i(g)=rd_i(g)+c_i.\end{equation} 

With these definitions we recall a set of monomials in $\mathbb{C}[\bold{x}_n]$ that descend to a vector-space basis for $S_n$ provided by Bango and Biagoli\cite{BB}.

	\begin{dfn}

		Given a r-colored permutation $g=\pi_1^{c_1}\pi_2^{c_2}\ldots \pi_n^{c_n}$, we define the {\bf r-descent monomial} $b_g$ as follows:\begin{equation}b_g:=\prod_{i=1}^n x^{f_i(g)}_{\pi_i}\end{equation}

	\end{dfn}

	The set of $r$-descent monomials descend to a basis for $S_n$.

We note that by the definition of $f_i$ and descents of r-colored permutations that $f_i(g)$ is a weakly decreasing sequence such that $f_i(g)-f_{i+1}(g)\leq r$.  

Chan and Rhoades\cite{CR} generalized these monomials to a set of monomials that descends to a basis for $S_{n,k}$.

	\begin{dfn}

		Given an r-colored permutation $g=\pi_1^{c_1}\pi_2^{c_2}\ldots \pi_n^{c_n}$ such that $des(g)<k$, and an integer sequence $I=(i_1, \ldots i_{n-k})$ such that $k-des(g)>i_1\geq i_2\geq\ldots\geq  i_{n-k} \geq 0$, we define the {\bf $(n,k,r)$-descent monomial} as\begin{equation}b_{g,I}:= b_g\cdot x^{ri_1}_{\pi_1}\ldots x^{ri_{n-k}}_{\pi_{n-k}}\end{equation}

	\end{dfn} 

	The set of $(n,k,r)$-descent monomials descend to a basis for $S_{n,k}$.
	We note that these monomials have individual powers strictly bounded by $kr$.  These observations motivate the following definitions. 

\begin{dfn}

	We call a partition an {\bf $(n,k,r)$-partition} if it has $n$ parts (some of which might be zero), each of which is strictly less than $rk$.

\end{dfn}

\begin{dfn}
	Given a partition $\mu$, we call an index $i$ an {\bf $r$-descent} of $\mu$ if\begin{equation}\lfloor \frac{\mu_i}{r}\rfloor>\lfloor \frac{\mu_{i+1}}{r}\rfloor.\end{equation}We will denote $Des^r(\mu)$ as the set of r-descents of $\mu$.

\end{dfn}

\begin{dfn}
	We call a partition an  {\bf r-descent partition} if the difference between consecutive parts is at most r, and the last part has size less than $r$.  
\end{dfn}

	The exponent partitions of both $r$ and $(n,k,r)$-descent monomials are $(n,k,r)$ partitions.  Furthermore the exponent partition of an $r$-descent monomials is an r-descent partition with r-descents equal to $Des(g)$ for the corresponding r-colored permutation $g$.  The last $k$ parts of the exponent partition of an $(n,k,r)$ descent monomial is an r-descent partition with r-descents determined by $Des_{n-k+1,n}(g)$, and the first $k$ parts have r-descents that are a superset of $Des_{1,n-k}(g)$.  Furthermore it is straightforward to see that all such $(n,k,r)$ partitions arise as the exponent partition of some $(n,k,r)$ descent monomial.

	In order to work with the basis of $(n,k,r)$-descent monomials we need to relate them to the partial order on monomials from the previous section. To do that we first give a way of associating an r-colored permutation to a monomial.

	\begin{dfn}

		Given a monomial $x_1^{a_1}\ldots x_n^{a_n}$ we define its {\bf index r-colored permutation} $g(m)=\pi_1^{c_1}\pi_2^{c_2}\ldots \pi_n^{c_n}$ to be the unique r-colored permutation such that

		\begin{enumerate}

			\item $a_{\pi_i}\geq a_{\pi_{i+1}}$ for $1\leq i<n$

			\item if $a_{\pi_i}=a_{\pi_{i+1}}$, then $\pi_i<\pi_{i+1}$

			\item $a_i\equiv c_i ($mod $ r)$

		\end{enumerate}

	\end{dfn}

	One last definition before we state some results is the following:

	\begin{dfn}

		Given a monomial $m=x_1^{a_1}\ldots x_n^{a_n}$ the {\bf r-complementary} partition $\mu(m)$ is the partition conjugate to\begin{equation}(\frac{a_{\pi_1}-rd_1(g)-c_1(g)}{r},\ldots,\frac{a_{\pi_n}-rd_n(g)-c_n(g)}{r})\end{equation}where $g=g(m)$ and $\pi$ is uncolored permutation of $g$. 

	\end{dfn}

	Implicit in these definitions is that they are well defined which is covered in \cite{BB}.

	We can now state the lemma that ties these objects together.

\begin{lem}

\label{Order2}

Let $m$ be a monomial equal to $x_1^{p_1}\ldots x_n^{p_n}$, then
among the monomials appearing in  $m\cdot e_{\mu}(\bold{x}_n^r)$,  the monomial\begin{equation}\prod_{i=1}^n x_{\pi(i)}^{p_{(\pi(i))}+ru'_i}\end{equation}is the maximum with respect to $\prec$, where $\pi$ is the index permutation of $m$.
	\begin{proof}
	The proof of this is similar to the proof of Prop \ref{Order}.
	\end{proof}

\end{lem}

Setting aside this side of things for a moment, we turn towards the representation side.
The analog of partitions which index the irreducible representations of $\mathfrak{S}_n=G(1,1,n)$ are r-partitions.

	\begin{dfn}
		An {\bf r-partition} of $n$ is an r-tuple of partitions $(\mu^0,\mu^1,\mu^2,\ldots, \mu^{r-1})$ such that $\sum_{i=0}^{r-1} \vert \mu_i\vert=n$.
		We will use Greek letters with a bar to denote r-partitions, and will write $\overline{\mu}\vdash_r n$ to denote that $\overline{\mu}$ is an r-partition of $n$.

	\end{dfn}

The conjugacy classes of $G_n$, and thus the irreducible representations of $G_n$, are indexed by r-partitions of $n$.  Given an r-partition $\overline{\lambda}$, we denote the irreducible representation of $G_n$ corresponding to $\overline{\lambda}$ as $S^{\overline{\lambda}}$.  The analog of standard Young tableaux in $G_n$ are standard Young r-tableaux.

	\begin{dfn}
		A {\bf standard Young r-tableau} of shape $(\mu^0,\mu^1,\ldots, \mu^{r-1})=\overline{\mu}$ is a way of assigning the integers $1,2,\ldots, n$ to the boxes of $r$ Ferrers diagrams of shapes $\mu^0, \mu^1,\ldots \mu^{r-1}$ such that in each of the Ferrers diagrams the integers increase down columns and along rows.  We denote the set of all standard Young r-tableaux of shape $\overline{\mu}$ as $SYT(\overline{\mu})$

	\end{dfn}

	As with standard Young tableaux we have the notion of descents.

	\begin{dfn}

		An index $i$ is a {\bf descent} of a standard Young r-tableaux $T$ if one of the following holds:

		\begin{enumerate}
			
			\item $i+1$ is in a component with a higher index than $i$

			\item $i+1$ and $i$ are in the same component and $i+1$ is strictly below $i$.

		\end{enumerate}
		
		Similar to other descents, we will denote $Des(T)$ as the set of all descents of $T$, $des(T)$ will be the number of descent of $T$, $d_i(T)$ will be the number of descents of $T$ that are $i$ or bigger.  One last statistic related to descents is\begin{equation}f_i(T)=r\cdot d_i(T)+c_i(T)\end{equation} where $c_i(T)$ is the index of the component of $T$ that contains $i$. 
	\end{dfn}

The result connecting standard Young r-tableaux to our problem is the following:

	\begin{prop}
	\label{BBProp}
		The graded trace of the action of $\overline{\tau}$ on $S_n$ has the following formula 
		\begin{equation}Tr_{\mathbb{C}[\bold{x}_n]}(\tau)=\frac{1}{\prod_{i=1}^n(1-q_1^rq_2^r\ldots q_i^r)}\sum_{\overline{\lambda}\vdash_r n}\chi_{\overline{\tau}}^{\overline{\lambda}}\sum_{T\in SYT(\overline{\lambda})}\prod_{i=1}^{n} q_i^{f_i(T)}\end{equation}

		Where $\chi_{\overline{\tau}}^{\overline{\lambda}}$ is the character $S^{\overline{\lambda}}$ evaluated at an element of type $\overline{\tau}$.

	\end{prop}

The proposition is proved in \cite{BB}, though it is in a more general form since the formula that they give is for the entire family of groups $G(r,p,n)$.  The formula that we state here is how it simplifies in the case $p=1$.

\subsection{Results}

	In order to calculate the multiplicities  of $S^{\overline{\lambda}}$ in $S_{n,k,\mu}$ we will calculate the graded trace of the action on an element of type $\overline{\tau}$ on the space of polynomials in $\mathbb{C}[\bold{x}_n]$ where the individual exponents of each variable are less than $kr$.  We denote this space $\mathbb{C}_{kr}[\bold{x}_n]$.  The first way we will do this calculaltion will be by using Proposition \ref{BBProp}.  The other way will use a basis for $\mathbb{C}_{kr}[\bold{x}_n]$ created from the descnent basis for $S_{n,k}$ .

\begin{lem}

\label{OrderLem2}
If  $m=\prod_{i=1}^n x_i^{a_i}$ is a monomial in $\mathbb{C}_{rk}[\bold{x}_n]$ (that is $a_i<kr$ for all $i$), then\begin{equation}m=b_{g,I}e_\nu(\bold{x}_n^r)+\sum.\end{equation} Where $g=g(m)$; $\sum$ is a sum of monomials $m'\prec m$; $I$ is a sequence defined by $i_\ell=\mu'_{\ell}-\mu'_{n-k+1}$ where $\mu=\mu(g)$; and $\nu$ is the partition specified by:

	\begin{enumerate}

		\item $\nu'_\ell=\mu'_\ell$ for $\ell> n-k$

		\item $\nu'_\ell=\mu'_{n-k+1}$ for $\ell\leq n-k$
	
	\end{enumerate}

Furthermore $\nu$ consists of parts of size at least $n-k+1$.

	\begin{proof}
	
	In order to show that $b_{g,I}$ is well defined we need to show that $k-des(g)>i_1\geq \ldots\geq i_{n-k}\geq 0$.  Since $I$ is formed by taking a weakly descreasing, non-negative sequence and subtracting a constant that is smaller than the smallest part, I satisfies $i_1\geq i_2\geq\ldots \geq i_{n-k}\geq 0$.  Letting $\pi$ be the uncolored permutation of $g$ and using definitions we get\begin{equation}i_1=\mu'_1-\mu'_{n-k+1}\leq \mu'_1=\frac{a_{\pi(1)}-rd_1(g)-c_1(g)}{r}\end{equation}and then by assumption $a_{\pi(1)}<rk$, and by definition $d_1(g)=des(g)$, thus\begin{equation}i_1< \frac{rk-rdes(g)}{r}=k-des(g).\end{equation} Therefore the use of $b_{g,I}$ is well defined.

	Next we show that $b_{g,I}$ and $m$ have the same index r-colored permutation, that is that\begin{equation}g(b_{g,I})=g(m)=g.\end{equation} To show this we consider the sequence of the exponents of $x_{\pi(\ell)}$ in $b_{g,I}$.  This sequence is the sum of the sequences $rd_\ell(g)+c_\ell(g)$ and $ri_\ell$ (where we take $i_\ell=0$ for $\ell>n-k$).  Since these are both weakly-decreasing sequences, their sum is also weakly-decreasing.  Furthermore if the $\ell$th and $(\ell+1)$th entries are equal then $rd_\ell(g)+c_\ell(g)=rd_{\ell+1}(g)+c_{\ell+1}(g)$, which means that $d_\ell(g)=d_{\ell+1}(g)$ and $c_\ell(g)=c_{\ell+1}(g)$ which by the definition of $d_\ell(g)$ implies that $\ell$ is not a descent of $g$.  Since $c_\ell(g)=c_{\ell+1}(g)$ this means that $\pi(\ell)<\pi(\ell+1)$, thus $g$ satisfies the first two conditions of being the index r-colored permutation, and the 3rd condition is immediate from the definition.

	Now by Proposition \ref{Order2}, the maximum monomial in $b_{g,I}e_\nu(\bold{x}_n^r)$ will have the form $\prod_{\ell=1}^n x_{\pi(\ell)}^{q_{\ell}}$ where $q_\ell$ is given by:

	\begin{enumerate}

		\item $q_\ell=rd_\ell(g)+c_\ell(g)+ri_\ell+r\nu'_\ell$ for $\ell\leq n-k$

		\item $q_\ell=d_\ell(g)+c_\ell(g)+r\nu'_\ell$ for $\ell> n-k$

	\end{enumerate}

By substitution, first using the definitions of $i_\ell$ and $\nu_\ell$ and then the definition of the r-complementary partition, we get that\begin{equation}q_\ell=d_\ell(g)+c_\ell(g)+r\mu'_\ell-r\mu'_{n-k+1}+r\mu'_{n-k+1}=d_\ell(g)+c_\ell(g)+r\mu'_\ell=a_{\pi(\ell)}\end{equation}for $\ell\leq n-k$, and\begin{equation}q_\ell=d_\ell(g)+c_\ell(g)-r\mu'_\ell=a_{\pi(\ell)}\end{equation}
for $\ell>n-k$

Finally $\nu$ has parts of size at least $n-k+1$ because by definition, the first $n-k+1$ parts of $\nu'$ are all the same size.

	\end{proof}
\end{lem}

\begin{prop}
\label{BasisProp2}

	The set $B_{n,k}$ consisting of products $b_{g,I}e_\nu(\bold{x}_n^r)$ where $\nu$ is a partition with parts of size at least $n-k+1$ and $(\lambda(b_{g,I})+r\nu')_1<rk$ form a basis for $\mathbb{C}_{rk}[\bold{x}_n]$.
	\begin{proof}

The condition that $(\lambda(b_{g,I})+r\nu')_1<rk$ along with \ref{OrderLem2} guarantees that each of the elements of $B_{n,k}$ are in $\mathbb{C}_{nk}[\bold{x}_n]$.

Iteratively applying Lemma \ref{OrderLem2} lets us express any monomial in $\mathbb{C}_{kr}[\bold{x}_n]$ as a linear combination elements of $B_{n,k}$, thus $B_{n,k}$ spans $\mathbb{C}_{rk}[\bold{x}_n]$.    To show that this expansion is unique (up to rearrangement) it is sufficient to show that if the maximal monomials in $b_{g,I}e_{\nu}(\bold{x}_n^r)$ and $b_{h,J}e_{\rho}(\bold{x}_n^r)$ are the same, then $g=h$, $I=J$ and $\nu=\rho$.  To see this, we note that as a corollary of the proof of Lemma \ref{OrderLem2}, the index r-colored permutations of the maximal monomials are the same, and they are both $g$ and $h$, and thus $g=h$.  Then, by Proposition \ref{Order2}, the power of $x_{\pi(\ell)}$ in each of these maximum monomials will be $rd_\ell(g)+c_\ell(g)+ri_\ell+r\nu'_\ell$ and $rd_\ell(h)+c_\ell(h)+rj_\ell+r\rho'_\ell$.  This immediately gives that $\nu'_\ell=\rho'_\ell$ for $\ell>n-k$ since $i_\ell=j_\ell=0$ for $\ell>n-k$.  Then since the first $n-k+1$ parts of $\nu'$  are all equal and the first $n-k+1$ parts of $\rho'$ are equal and since $\nu'_{n-k+1}=\rho'_{n-k+1}$, we have that $\nu'=\rho'$ which implies $\nu=\rho$.  This then implies that $i_\ell=j_\ell$ for all $\ell$, and therefore this expansion is unique.  Therefore $B_{n,k}$ is linearly independent and is a basis.

	\end{proof}

\end{prop}

\begin{prop}
\label{DomProp2}
	Let $p$ be the map projecting from $\mathbb{C}[\bold{x}_n]$ to $S_{n,k}$ and let $m$ be a monomial in $\mathbb{C}_{rk}[\bold{x}_n]$.  Then\begin{equation}p(m)=\sum_{g,I} \alpha_{g,I} b_{g,I}\end{equation}where $\alpha_{g,I}$ are some constants, and the sum is over pairs $g,I$ such that $\lambda(b_{g,I})\trianglelefteq \lambda(m)$.

	\begin{proof}
		Since $B_{n,k}$ is a basis we can express $m=\sum_{g,I,\nu} \alpha_{g,I,\nu} b_{g,I} e_\nu(\bold{x}_n^r)$. By Lemma \ref{OrderLem2} $\alpha_{g,I,\nu}$ is zero if the leading monomial of $b_{g,I}e_\nu$ is not weakly smaller than $m$ under the partial order on monomials.  But since the partial order on monomials refines the dominance order on exponent partitions, for each non-zero term the exponent partition of the leading monomial will be dominated by $\lambda(m)$ that is that\begin{equation}(\lambda(b_{g,I})+r\nu')\trianglelefteq \lambda(m).\end{equation} Then when we project down to $S_{n,k}$, each term with $\nu\neq \emptyset$ will vanish since $e_\nu(\bold{x}_n^r)$ is in $J_{n,k}$, so that\begin{equation}p(m)=\sum_{g,I}\alpha_{g,I,\emptyset} b_{g,I}\end{equation}where the sum is over $(g,I)$ such that $\lambda(b_{g,I})\trianglelefteq \lambda(m)$.
	\end{proof}

\end{prop}

Proposition \ref{DomProp2} give the following Corollary:

\begin{coro}
\label{zero2}

	$S_{n,k,\rho}$ is zero unless $\rho$ is the exponent partition of an $(n,k,r)$-descent monomial, which occurs precisely when $\rho$ is an $(n,k,r)$ partition such that the last $k$ parts form an $r$-descent partition.

\end{coro}

This basis allows us to express the action of $\tau\in G_n$ on $ \mathbb{C}_{rk}[\bold{x}_n]$ in terms of its action on $S_{n,k}$ with the basis of $(n,k,r)$-descent monomials.  Specifically if\begin{equation}\tau(b_{g,I})=\sum_{h,J} \alpha_{h,J} b_{h,J},\end{equation}for some constants $\alpha_{h,J}$, then\begin{equation}\tau(b_{g,I} e_\nu(\bold{x}_n^r))=\sum_{h,J} \alpha_{h,J} b_{h,J} e_\nu(\bold{x}_n^r).\end{equation}

This equality holds, since $e_\nu(\bold{x}_n^r)$ is invariant under the action of $G_n$.  The important part is that for each $(h,J)$ such that $\alpha_{h,J}$ is non-zero, $b_{h,J}e_\nu(\bold{x}_n^r)$ is an element of our basis, that is we do not violate the condition that $(\lambda(b_{h,J})+r\nu')_1<rk$.  By Proposition \ref{DomProp2} we have that $\lambda(b_{h,J})\triangleleft \lambda(b_{g,I})$ which implies that $\lambda(b_{h,J})+r\nu'\triangleleft \lambda(b_{g,I})+r\nu'$, and thus $(\lambda(b_{h,J})+r\nu')_1\leq (\lambda(b_{g,I})+r\nu')_1<rk$, therefore it satisfies the condition to be in our basis.

\begin{lem}

\label{DecomLem2}
 Given an $(n,k,r)$-partition $\mu$ and an $(n,k,r)$ r-descent partition $\nu$ there exists a unique $(n,k,r)$ partition $\rho$ such that $\mu=\nu+r\rho$ if and only if $Des^r(\nu)\subseteq Des^r(\mu)$ and $\mu_i\equiv \nu_i$ (mod r) for all $i$.

	\begin{proof}
		There is only one possible value for each part of $\rho$ which is $\rho_i=\frac{\mu_i-\nu_i}{r}$.  The mod $r$ condition is necessary and sufficient for these values to be integers.  In order for this to be a partition we need\begin{equation}\rho_i-\rho_{i+1}=\frac{1}{r}[(\mu_i-\mu_{i+1})-(\nu_i-\nu_{i+1})]\geq 0.\end{equation} Let $c_i$ be the common remainder of $\mu_i$ and $\nu_i$ mod r. Since $\nu$ is an  r-descent partition, $\frac{1}{r}((\nu_i-c_i)-(\nu_{i+1}-c_{i+1}))$ is 1 if $i$ is an r-descent of $\nu$ and 0 if it is not.  Similarly, $\frac{1}{r}((\mu_i-c_i)-(\mu_{i+1}-c_{i+1}))$ is at least 1 if $i$ is an r-descent of $\mu$ and 0 otherwise.  Thus in order for\begin{equation}\frac{1}{r}[(\mu_i-c_i)-(\mu_{i+1}-c_{i+1})]-\frac{1}{r}[(\nu_i-c_i)-(\nu_{i+1}-c_{i+1})]=\frac{1}{r}[(\mu_i-\mu_{i+1})-(\nu_i-\nu_{i+1})]\end{equation}to be non-negative, it is necessary and sufficient that if $i$ is an r-descent of $\nu$, then $i$ is also an r-descent of $\mu$.  That is, $\rho$ will be a partition if and only if $Des^r(\nu)\subseteq Des^r(\mu)$.
	\end{proof}
\end{lem}

\begin{lem}

\label{Decom2Lem2}

Given an $(n,k,r)$ partition $\mu$ and a set $S\subseteq Des^r_{n-k+1,n}(\mu)$, there is a unique pair $(\nu,\rho)$ such that $\mu=\nu+r\rho$ and $\nu$ is the exponent partition of an $(n,k,r)$-descent monomial with $Des^r_{n-k+1,n}(\nu)=S$, and $\rho$ is an $(n,k)$-partition with $\rho_1=\rho_2=\ldots =\rho_{n-k+1}$.

	\begin{proof}

		The last $k$ values of the exponent partition of an $(n,k,r)$ descent monomial form an r-descent partition, so applying Lemma \ref{DecomLem2} to the partition determined by $S$ and the values of $\mu_i$ mod(r) determines the last $k$ values of $\rho$.  But since we need that the first $n-k+1$ values of $\rho$ are the same, this determines what $\rho$ must be, and thus by subtraction what $\nu$ must be.  We just need to check that $\nu$ is actually a partition, that is that $\nu_i-\nu_{i+1}\geq 0$ for $1\leq i\leq n-k$.  This is true since $\nu_i-\nu_{i+1}=(\mu_{i}-r\rho_i)-(\mu_{i+1}-r\rho_{i+1})=\mu_i-\mu_{i+1}\geq 0$ since $\rho_i=\rho_{i+1}$.

	\end{proof}

\end{lem}

We now give the proof of Theorem \ref{Result2}

	\begin{proof}(Thm \ref{Result2})
		The condition on when $S_{n,k,\rho}$ is zero is covered by Corollary \ref{zero2}.

		 We consider the graded trace of the action of $\tau\in G_n$ on $\mathbb{C}_{rk}[\bold{x}_n]$ defined by\begin{equation}Tr_{\mathbb{C}_{rk}[\bold{x}_n]}(\overline{\tau}):=\sum_m \langle \overline{\tau}(m),m\rangle \cdot \bar{q}^{\lambda(m)}.\end{equation}  From \ref{BBProp} we have that\begin{equation}Tr_{\mathbb{C}[\bold{x}_n]}(\tau)=\frac{1}{\prod_{i=1}^n(1-q_1^rq_2^r\ldots q_i^r)}\sum_{\overline{\lambda}}\chi_{\overline{\mu}}^{\overline{\lambda}}\sum_{T\in SYT(\overline{\lambda})}\prod_{i=1}^{n} q_i^{f_i(T)}\end{equation} (where $\overline{\mu}$ is the cycle type of $\tau$).
From this we can recover $Tr_{\mathbb{C}_{rk}[\bold{x}_n]}(\tau)$ by restricting to powers of $q_1$ that are at most $rk-1$.  Doing this gives\begin{equation}\sum_{\overline{\lambda}\vdash_r n}\chi_{\overline{\mu}}^{\overline{\lambda}}\sum_{T\in SYT(\overline{\lambda}),\nu}\bar{q}^{\lambda_{F(T)}}\bar{q}^{r\nu}.\end{equation} Where $F(T)=(f_1(T),f_2(T),\ldots, f_n(T))$ and the $\nu$'s are partitions such that $(\lambda_{F(T)})_1+r\nu_1<rk$.

		Alternatively, we can calculate $Tr_{\mathbb{C}_{rk}[\bold{x}_n]}(\tau)$ by using the basis from Proposition \ref{BasisProp2}, this gives
\begin{align}Tr_{\mathbb{C}_{rk}[\bold{x}_n]}(\tau)&=\sum_{g,I,\nu} \langle \tau(b_{g,I} e_\nu(\bold{x}_n^r),b_{g,I} e_\nu(\bold{x}_n^r) \rangle \bar{q}^{\lambda(b_{g,I})}\bar{q}^{r\nu'}\\
&=\sum_{g,I,\nu} \langle \tau(b_{g,I}),b_{g,I} \rangle \bar{q}^{\lambda(b_{g,I})}\bar{q}^{r\nu'}\\
&=\sum_{\phi,\nu} Tr_{S_{n,k}}(\tau; \bar{q}^\phi)\bar{q}^{\phi}\bar{q}^{r\nu'}\end{align}
where the $\nu$'s are partitions with parts of size at least $n-k+1$ such that $(\phi)_1+(r\nu')_1<rk$, the $\phi$'s are the exponent partitions of $(n,k,r)$-descent monomials, and $Tr_{S_{n,k}}(\tau;\bar{q}^\phi)$ is the coefficient of $\bar{q}^\phi$ in the graded trace of the action of $\tau$ on $S_{n,k}$.

		We now consider the coefficient of $\bar{q}^\rho$ for some $(n,k,r)$-partition $\rho$.  Using the first calculation and Lemma \ref{DecomLem2}, the inner sum can be reduced to $T$ such that $Des(T)\subseteq Des^r(\rho)$, and such that $c_i(T)\equiv (\rho)_i$ mod $r$ for all $i$, so that we get\begin{equation}\sum_{\overline{\lambda}\vdash_r n} \chi^{\overline{\lambda}}_{\overline{\mu}}\vert \{T\in SYT(\lambda): Des(T)\subseteq Des^r(\rho)\text{, and } c_i(T)\equiv \rho_i (\text{mod } r)\}\vert.\end{equation}

Looking at the second calculation and using Lemma \ref{Decom2Lem2} gives\begin{equation}\sum_{\phi} Tr_{S_{n,k}}(\tau;\bar{q}^{\phi}),\end{equation}where the sum is over a set consisting of the exponent partitions $\phi$'s of $(n,k,r)$ descent monomials with the set of $\phi$'s being in bijection with subsets of $Des^r_{n-k+1,n}(\rho)$ such that for each $\phi$, $Des^r_{n-k+1,n}(\phi)$ is equal to the corresponding subset, and $\phi_i\equiv \rho_i$ mod r for all $i$.  Together this gives that\begin{equation}\sum_{\overline{\lambda}\vdash_r n} \chi^{\overline{\lambda}}_{\overline{\mu}}\vert \{T\in SYT(\overline{\lambda}): Des(T)\subseteq Des^r(\rho), c_i(T)\equiv \rho_i (\text{mod } r)\}\vert=\sum_{\phi} Tr_{S_{n,k}}(\tau;\bar{q}^{\phi}).\end{equation}

	We want to further refine this result by showing that\begin{equation}\sum_{\overline{\lambda}\vdash_r n} \chi^{\overline{\lambda}}_{\overline{\mu}}\vert \{T\in SYT(\overline{\lambda}): Des_{n-k+1,n}(\phi)\subseteq Des(T)\subseteq Des(\phi),c_i(T)\equiv \rho_i (\text{mod }r) \}\vert\end{equation} \begin{equation}= Tr_{S_{n,k}}(\tau;\bar{q}^{\phi})\end{equation}for any specific $\phi'$.  We do this by induction on $\vert \phi' \vert$.  The base case of $\phi'$ being empty is trivial.  If we take $\rho=\phi'$, then $\phi'$ will appear in the sum, and all other $\phi$'s will be smaller, so by the inductive hypothesis,\begin{equation}\sum_{\overline{\lambda}\vdash_r n} \chi^{\overline{\lambda}}_{\overline{\mu}}\vert \{T\in SYT(\overline{\lambda}): Des_{n-k+1,n}(\phi)\not\subseteq Des(T)\subseteq Des(\phi),c_i(T)\equiv \rho_i (\text{mod }r) \}\vert\end{equation} \begin{equation}= \sum_{\phi\neq \phi'}Tr_{S_{n,k}}(\tau;\bar{q}^{\phi})\end{equation}
Subtracting this from our result gives the desired refinement.  This then proves the Theorem since the exponent partition of any $(n,k,r)$-descent monomials $b_{g, I}$ appears when we take  $\rho=\lambda(b_{g,I})$.

	\end{proof}

As an example of Theorem \ref{Result2}, let $n=7$, $k=5$ and let $r=2$.  Then consider letting $\rho=(9,5,5,4,3,2,0)$.  The standard Young r-tableaux $T$ that we must consider will have $4,6,7$ in the 0-componet, and $1,2,3,5$ in the 1-component.  Furthermore they will have  $\{4,6\}\subseteq Des(T)\subseteq \{1,4,6\}$. The possibilities for the 0-component and the 1-component are independent.  The possibilities for the 0-component are

\medskip

	\ytableaushort{46,7}*[*(blue!20)]{2}\, ,
	\ytableaushort{4,6,7}*[*(blue!20)]{1,1}.\,
	\medskip

The possibilities for the 1-componet are

\medskip

	\ytableaushort{1235}\,,
	\ytableaushort{123,5}\,,
	\ytableaushort{135,2}*[*(blue!20)]{1}\,,
	\ytableaushort{13,25}*[*(blue!20)]{1}\,,	
	\ytableaushort{13,2,5}*[*(blue!20)]{1}\,.
	\medskip

Therefore the multiplicites of $S^{(2,1),(4)}$, $S^{(2,1),(2,2)}$, $S^{(2,1),(2,1,1)}$, $S^{(1,1,1),(4)}$,\\ $S^{(1,1,1),(2,2)},$ and $S^{(1,1,1),(2,1,1)}$ in $S_{n,k,\rho}$ are 1, and the multiplicities of $S^{(2,1),(3,1)}$ and $S^{(1,1,1),(3,1)}$ are 2. All other multiplicities are zero.

\medskip

Theorem \ref{Result2} also allows us to recover the following Corollary which is eqvuivalent (thier indexing is different) to a result to of Chan and Rhoades \cite{CR}.

\begin{coro}
	Let $f_{\overline{\lambda}}(q)$ be the generating function for the multiplicities of $S^{\overline{\lambda}}$ in the degree $d$ compent of $S_{n,k}$.  Then\begin{equation}f_{\overline{\lambda}}(q)=\sum_{T\in SYT(\overline{\lambda})} q^{maj(T)}\genfrac{[}{]}{0pt}{}{n-des(T)-1}{n-k}_{q^r}.\end{equation}
	Where the major index $maj(T)$ is equal to $\sum_{i=1}^n rd_i(T)+c_i(T)$
	\begin{proof}

		By Theorem \ref{Result2}, each standard Young r-tableau of shape $\overline{\lambda}$ contributes to $f_{\overline{\lambda}}(q)$ once for each partition $\rho$ such that $\rho$ is the exponent partition of an $(n,k,r)$-descent monomial and $Des^r_{n-k+1,n}(\rho)\subseteq Des(T)\subseteq Des^r(\rho)$.  All such $\rho$ come from $(n,k,r)$-descent monomials $b_{g,I}$ where $g$ is an $r$-colored permutation with $Des^r(g)=Des(T)$, $c_i(g)=c_i(T)$ for all $i$, and $I$ is a sequence such that $k-des(T)>i_1\geq i_2\geq\ldots\geq i_{n-k}\geq 0$.  This choice of $I$ is the same as choosing a partition that fits in an $(n-k)\times (k-1-des(T))$ box.  The generating function for the number of partitions of size $d$ that fit in an $(n-k)\times (k-1-des(T))$ box is $\genfrac{[}{]}{0pt}{}{(n-k)+(k-des(T)-1)}{n-k}_q=\genfrac{[}{]}{0pt}{}{n-des(T)-1}{n-k}_q$.  But in $b_{g,I}$ we are multiplying the values in $I$ by $r$, so we need to plug $q^r$ into this $q$-binomial coefficient to get $\genfrac{[}{]}{0pt}{}{n-des(T)-1}{n-k}_{q^r}$.   The factor of $b_g$ in the $(n,k,r)$ descent monomial then has degree $maj(T)$, so that each standard Young tableau $T$ of shape $\lambda$ will contribute $q^{maj(T)}\genfrac{[}{]}{0pt}{}{n-des(T)-1}{n-k}_{q^r}$  to $f_{\overline{\lambda}}(q)$.  This completes the proof.

	\end{proof}
\end{coro}

The proof of this result in \cite{CR} is fairly involved using a tricky recursive argument involving an auxillary family of algebras.  We manage to give a simpler proof for this result.

Overlapping notations slightly, Chan and Rhoades \cite{CR} also defined the ideal\begin{equation}I_{n,k}:=\langle x_1^{kr+1},x_2^{kr+1},\ldots, x_n^{kr+1}, e_n(\bold{x}_n^r),e_{n-1}(\bold{x}_n^r),\ldots, e_{n-k+1}(\bold{x}_n^r)\rangle,\end{equation}and the algebra\begin{equation}R_{n,k}:=\frac{\mathbb{C}[x_1,x_2,\ldots, x_n]}{I_{n,k}}.\end{equation} As before we can refine the grading on this algbra to define $R_{n,k,\rho}$, and can ask what the the graded isomorphism type of this $G_n$ module is.  By slightly modifying the results of this section (looking at partitions with largest part $kr$ instead of $kr-1$, and using the extended descent monomials from \cite{CR} instead of the descent monomials) we can a result that is analogous to Theorem \ref{Result2}.

\section{Conclusion}
\label{Conclusion}

One path to take from here would be to try to extend the $(n,k)$ coinvariant algebras introduced by Chan and Rhoades \cite{CR} for $G(r,1,n)$ to all complex reflection groups. It seems that the simplest groups to consider are $G(2,2,n)$ which are equal to the real reflection groups of Coxeter-Dynkin type $D_n$.  In the case that $G(r,1,n)$ is a real reflection group, the structure of the corresponding $(n,k)$ coinvariant algebra is governed by the combinatorics of the $k$-dimensional faces of the associated Coxeter complex.  We can define a candidate graded algebras for $D_n$ that will satisfy this property by using a more general technique of Garsia and Procesi \cite{GP} which we recall here	.

We start by taking a finite set of points $X\subset \mathbb{C}^n$.  We then consider the set of polynomials in $\mathbb{C}[x_1,x_2,\ldots, x_n]$ that vanish on $X$, that is\begin{equation}\{f\in \mathbb{C}[x_1,\ldots, x_n]: f(\overline{x})=0 \  \text{ for all } \ \overline{x}\in X\}.\end{equation}  This set is an ideal in $\mathbb{C}[x_1,\ldots, x_n]$, and we will denote it by $\mathbb{I}(X)$.  Next, we consider the quotient $\frac{\mathbb{C}[x_1,\ldots, x_n]}{\mathbb{I}(X)}$.  

The elements of this quotient can be viewed as $\mathbb{C}$-valued function on $X$.  We do this by taking a representative polynomial in $\mathbb{C}[x_1,\ldots, x_n]$, viewing it as a function from $\mathbb{C}^n$ to $\mathbb{C}$, and then restricting its domain to $X$.  Two polynomials will give rise to the same $\mathbb{C}$-valued function if and only if their difference vanishes on $X$, which occur precisely if the difference is in $\mathbb{I}(X)$.  Therefore this is well defined.  Furthermore for every element $\overline{x}\in X$, we can construct an indicator function for $\overline{x}$ as follows. For each element $\overline{y}\in X\backslash \overline{x}$ choose an index $i_{\overline{y}}$ at which $\overline{x}$ and $\overline{y}$ differ, then\begin{equation}\prod_{\overline{y}\in X\backslash \overline{x}} \frac{x_{i_{\overline{y}}}-\overline{y}_{i_{\overline{y}}}}{\overline{x}_{i_{\overline{y}}}-\overline{y}_{i_{\overline{y}}}}\end{equation}is an indicator function for $\overline{x}$. Therefore $\frac{\mathbb{C}[x_1,\ldots, x_n]}{\mathbb{I}(X)}$ is isomorphic as a vector space to $\mathbb{C}[X]^*\cong\mathbb{C}[X]$.

Any subgroup $W$ of $GL(\mathbb{C}^n)$ acts on $\mathbb{C}[x_1,\ldots, x_n]$ by linear substitution. If $X$ is invariant under $W$, then $\mathbb{I}(X)$ is invariant under $W$, and thus both$\frac{\mathbb{C}[x_1,\ldots, x_n]}{\mathbb{I}(X)}$ and $\mathbb{C}[X]$ are $W$-modules.  Furthermore in addition to being isomorphic as vector spaces, these two objects are isomorphics as $W$-modules.  Unfortunately, $\mathbb{I}(X)$ will not generally be homogeneous, and thus we will not have that $\frac{\mathbb{C}[x_1,\ldots, x_n]}{\mathbb{I}(X)}$ is graded.  In order to fix this we introduce a function $\tau$ that sends a non-zero polynomial to its top degree component.  For example\begin{equation}\tau(x_1^2+x_2^2+x_2x_3-x_1-x_2-x_3+3)=x_1^2+x_2^2+x_2x_3\end{equation}
and\begin{equation}\tau(x_1^4+x_1x_2x_3x_4+x_3^4-x_2^3-x_2^2+3)=x_1^4+x_1x_2x_3x_4+x_3^4.\end{equation}

We then consider the ideal $\mathbb{T}(X)$ generated by the top degrees of polynomials that vanish on $X$, that is\begin{equation}\mathbb{T}(X):=\langle\{\tau(f):f\in \mathbb{I}(X)-\{0\} \}\rangle.\end{equation} This ideal is homogeneous and invariant under $W$, therefore  $\frac{\mathbb{C}[x_1,\ldots, x_n]}{\mathbb{T}(X)}$ is a graded $W$-module.  Furthermore it can be shown that\begin{equation}\frac{\mathbb{C}[x_1,\ldots, x_n]}{\mathbb{T}(X)}\cong_W \frac{\mathbb{C}[x_1,\ldots, x_n]}{\mathbb{I}(X)}\cong_W \mathbb{C}[X].\end{equation} Then if we take $W$ to be $D_n$ and take $X$ to be a set of a points corresponding to the $k$-dimensional faces of the Coxeter complex of $D_n$ that is invariant under $D_n$,  our candidate algebra will be  $\frac{\mathbb{C}[x_1,\ldots, x_n]}{\mathbb{T}(X)}$.

There are two difficulties that we run into at this point. The first is the question of how we choose $X$.  Different choices of $X$ lead to isomorphic ungraded $D_n$-modules, but the graded structure in general depends on $X$, and it is not clear what the ``correct" choice is.  The second difficulty is getting a nice generating set for $\mathbb{T}(X)$.  We do have a general method to get a (potentially ugly) description of $\mathbb{T}(X)$ from $X$ which is the following.

The idea behind our method is that if we find a set  $P\subset \mathbb{T}(X)$ such that\begin{equation}dim(\frac{\mathbb{C}[x_1,\ldots, x_n]}{\langle P\rangle})=dim(\frac{\mathbb{C}[x_1,\ldots, x_n]}{\mathbb{T}(X)})=\vert X\vert,\end{equation}then $P$ generates $\mathbb{T}(X)$.  For a given $P\subset \mathbb{T}(X)$, let $st(P)$ be the standard monomial basis for$\frac{\mathbb{C}[x_1,\ldots, x_n]}{\langle P\rangle}$ with repsect to some graded monomial ordering (see \cite{Gr} for more details).  We will have found a $P$ that works when we have that $\vert st(P)\vert=\vert X\vert$. 

	All monomials of degree $\vert X\vert$ will appear in $\mathbb{T}(X)$.  In order to avoid cumbersome notation we will show this with an example.  If $n=3$ and $X$ consists of the points $(\alpha_1,\alpha_2,\alpha_3), (\beta_1,\beta_2,\beta_3),$ and $(\gamma_1,\gamma_2,\gamma_3)$, then have \begin{align}x_1^3&=\tau((x_1-\alpha_1)(x_1-\beta_1)(x_1-\gamma_1))\\ x_1^2x_2&=\tau((x_1-\alpha_1)(x_1-\beta_1)(x_2-\gamma_2))\\ x_1^2x_3&=\tau((x_1-\alpha_1)(x_1-\beta_1)(x_3-\gamma_3)),\end{align} and so on.
This idea generalizes to show that all degree $\vert X\vert$ monomials will appear in $\mathbb{T}(X)$.  We will thus start with $P$ consisting of all monomials of degree $d$.  Then $st(P)$ will consist of monomials of degree less than $d$, which is a finite set.

We now describe a method for adding an element to $P$ that will reduce the size of $st(P)$.  Let $m_1,m_2,\ldots m_s$ be the elements of $st(P)$, and let $p_1,p_2,\ldots p_t$ be the elements of $X$.  We then create a $t\times s$ matrix $M$ by setting $M_{ij}=m_j(p_i)$.  If $t<s$, then the null space of $M$ is non-zero, so we can take a non-zero vector $v=(v_1,v_2,\ldots, v_s)$ in the null space.  We then consider the polynomial $f=\sum_{j=1}^s v_jm_j$.  Evaluating this polynomial at $p_i$ gives $\sum_{j=1}^s v_jm_j(p_i)=(Mv)_i=0$.  Thus $f$ vanishes on $X$ which means that $\tau (f)$ is in $\mathbb{T}(X)$.  The leading monomial of $\tau (f)$ is an element of $st(P)$, and adding $\tau (f)$ to $P$ will at least eliminate this leading monomial from $st(P)$.  We then iterate this process until $\vert st(P)\vert =\vert X\vert$.

This method also gives us the standard monomial basis for $\frac{\mathbb{C}[x_1,\ldots, x_n]}{\mathbb{T}(X)}$.  This allows us to give examples of when different choices of $X$ lead to different graded structures.  If we let $X$ be the orbits of $(1,1,2)$, $(-1,1,2)$, and $(1,2,2)$ under $D_3$, then the Hilbert series of $\frac{\mathbb{C}[x_1,\ldots, x_n]}{\mathbb{T}(X)}$
is $5q^5+11q^4+10q^3+6q^2+3+1$.  If instead we take the orbits of $(1,1,2)$, $(-1,1,2)$, and $(1,\sqrt{\frac{5}{2}},\sqrt{\frac{5}{2}})$, then we get Hilbert series $11q^5+9q^4+7q^3+5q^2+3q+1$, and if we take the orbits of $(1,1,2)$, $(-1,1,2)$, and $(0,\sqrt{3},\sqrt{3})$ we get Hilbert series $4q^6+8q^5+8q^4+7q^3+5q^2+3q+1$.  From experimental data is does appear that there is a generic isomorphism type, but even for $X$ that give rise to isomorphic graded $D_n$-modules, the ideals $\mathbb{T}(X)$ will be different.

\section{Acknowledgements}  I would like to thank my advisor Brendon Rhoades for his work in directing and advising this work.

\bibliographystyle{plain}
\bibliography{MeyerDescentPrePrint2}

\begin{thebibliography}{10}

\bibitem{Gr}
William~Wells Adams and Philippe Loustaunau.
\newblock {\em An introduction to Gr{\"o}bner bases}.
\newblock Number~3. American Mathematical Soc., 1994.

\bibitem{ABR}
Ron Adin, Francesco Brenti, and Yuval Roichman.
\newblock Descent representations and multivariate statistics.
\newblock {\em Transactions of the American Mathematical Society},
  357(8):3051--3082, 2005.

\bibitem{BB}
Eli Bagno and Riccardo Biagioli.
\newblock Colored-descent representations of complex reflection groups g (r, p,
  n).
\newblock {\em Israel Journal of Mathematics}, 160(1):317--347, 2007.

\bibitem{Sh}
Georgia Benkart, Laura Colmenarejo, Pamela~E Harris, Rosa Orellana, Greta
  Panova, Anne Schilling, and Martha Yip.
\newblock A minimaj-preserving crystal on ordered multiset partitions.
\newblock {\em Advances in Applied Mathematics}, 95:96--115, 2018.

\bibitem{CR}
Kin Tung~Jonathan Chan and Brendon Rhoades.
\newblock Generalized coinvariant algebras for wreath products.
\newblock {\em arXiv preprint arXiv:1701.06256}, 2017.

\bibitem{C}
Claude Chevalley.
\newblock Invariants of finite groups generated by reflections.
\newblock {\em American Journal of Mathematics}, 77(4):778--782, 1955.

\bibitem{G}
Adriano~M Garsia.
\newblock Combinatorial methods in the theory of {Cohen-Macaulay} rings.
\newblock {\em Advances in Mathematics}, 38(3):229--266, 1980.

\bibitem{GP}
Adriano~M Garsia and Claudio Procesi.
\newblock On certain graded {Sn}-modules and the q-{Kostka} polynomials.
\newblock {\em Advances in Mathematics}, 94(1):82--138, 1992.

\bibitem{GS}
Adriano~M Garsia and Dennis Stanton.
\newblock Group actions on stanley-reisner rings and invariants of permutation
  groups.
\newblock {\em Advances in Mathematics}, 51(2):107--201, 1984.

\bibitem{HRS}
James Haglund, Brendon Rhoades, and Mark Shimozono.
\newblock Ordered set partitions, generalized coinvariant algebras, and the
  {Delta Conjecture}.
\newblock {\em arXiv preprint arXiv:1609.07575}, 2016.

\bibitem{S}
Richard~P Stanley.
\newblock Invariants of finite groups and their applications to combinatorics.
\newblock {\em Bulletin of the American Mathematical Society}, 1(3):475--511,
  1979.

\end{thebibliography}

\end{document}